\def\ZZ         {{\mathbb Z}}
\def\RR         {{\mathbb R}}
\def\CC         {{\mathbb C}}
\def\QQ         {{\mathbb Q}}
\def\PP         {{\mathbb P}}
\def\ii         {{\rm i}}
\def\ee         {{\rm e}}

\def\rk         {{\rm rk}}
\def\Spec         {{\rm Spec}}

\def\dim        {{\rm dim}}
\def\sin        {{\rm sin}}

\def\Ell        {{\cal ELL}}

\def\Prod       {{\displaystyle\prod}}

\def\cal        {\mathcal}

\documentclass{amsart}
\usepackage{amssymb}
\usepackage{verbatim}
\usepackage{eucal}
\usepackage{mathrsfs}
\usepackage{graphicx}
\usepackage{psfrag}

\newtheorem{theorem}{Theorem}[section]
\newtheorem{lemma}[theorem]{Lemma}
\newtheorem{proposition}[theorem]{Proposition}
\newtheorem{corollary}[theorem]{Corollary}

\theoremstyle{definition}
\newtheorem{definition}[theorem]{Definition}

\newtheorem{example}[theorem]{Example}

\theoremstyle{remark}
\newtheorem{remark}[theorem]{Remark}

\begin{document}

\title{McKay correspondence for elliptic genera}

\author{Lev Borisov}
\address{Department of Mathematics\\
Columbia University\\
New York, NY  10027 and }
\address{Department of Mathematics\\
University of Wisconsin-Madison}
\email{lborisov@math.columbia.edu}

\author{Anatoly Libgober}
\address{Department of Mathematics\\
University of Illinois\\
Chicago, IL 60607}
\email{libgober@math.uic.edu}

\thanks{The first author was partially supported by NSF grant DMS-0140172.
The second author was partially supoorted by NSF grant DMS-9803623}

\begin{abstract}
We establish a correspondence between orbifold and singular
elliptic genera of a global quotient. While the former is defined in
terms of the fixed point set of the action, the latter is defined in terms
of the resolution of singularities. As a byproduct, the second quantization
formula of Dijkgraaf, Moore, Verlinde and Verlinde is extended to 
arbitrary Kawamata log-terminal pairs.
\end{abstract}

\maketitle

\section{Introduction}

McKay correspondence was originally proposed in \cite{McKay} as a relation 
between minimal resolutions of quotient singularities ${\CC}^2/G$ 
where $G$ is a finite subgroup of $SL_2({\CC})$ and 
the representations of $G$. 
Shortly after that, Dixon, Harvey, Vafa and Witten (cf. \cite{DHVW}) 
discovered
a formula for the euler characteristic of certain resolutions of 
quotients:  

\begin{equation} e(\widetilde {X/G})=
\frac 1{\vert G \vert}\sum_{gh=hg} e(X^{g,h})
\label{eulercommuting}
\end{equation} 
where $X$ is a complex manifold, $\pi^*: \widetilde {X/G} \rightarrow X/G$ 
is a resolution of 
singularities such that 
$\pi^*K_{X/G}=K_{\widetilde {X/G}}$, $X^{g,f}$ is the submanifold
of $X$ of points fixed by both $f$ and $g$. The right hand side 
in (\ref{eulercommuting}) can be written as  the 
sum ofver the conjugacy classes: $\sum_{\{g\}} e(X^g/C(g))$
where $C(g)$ is the centralizer of $g$ which 
for $X={\bf C}^2$ 
is the number of irreducible representations of $G$. At the same time, 
the other side   
in (\ref{eulercommuting}) is 
the number of exceptional 
curves in minimal resolution plus $1$ and one obtains McKay 
correspondence on the numerology level (cf. \cite{HirzHofer}).  
McKay correspondence became the subject of intense study
and the term is now primarily used to indicate a relationship between the 
various invariants of the actions of finite automorphism groups on 
quasiprojective varieties and resolutions of the corresponding quotients
by such actions generalizing (\ref{eulercommuting}). We refer to the report 
\cite{Reid} for a survey of evolution of 
ideas since original empirical observation of McKay. 

One of the main results to date on the relationship between the invariants 
of actions and resolutions of quotients is the description of 
the $E$-function of 
a crepant resolution in terms of the invariants of the action 
(cf \cite{Batyrev}, \cite{denefl})). 
We recall that for a quasiprojective 
variety $M$ its  $E$-function is defined as 
$$E(M;u,v)\colon = \sum_{o,q} u^pv^q\sum_n (-1)^nh^{p,q}(H^n_c(M))$$
where $h^{p,q}(H^n_c(M))$ are the Hodge numbers 
of the Deligne's mixed Hodge structure on the compactly supported
cohomology of $M$.
$E$-function incorporates many classical numerical invariants of 
manifolds. For example, if $M$ is a projective manifold and $(u,v)=(y,1)$
one obtains Hirzebruch's $\chi_y$-genus which in turn has 
topological and holomorphic euler characteristics and signature as 
its special values. 

In \cite{Batyrev}, Batyrev extended the definition of 
$E$-function to the case of global quotient of a variety $M$ by 
a finite group $G$. 
He defined the orbifold $E$-function, 
$E_{orb}(M,G;u,v)$ in terms of the action of a finite group $G$.
Moreover, he extended this definition to $G$-normal pairs 
$(M,D)$ that consist of a variety $M$ and a simple normal 
crossing $G$-equivariant divisor $D$ on it. 
Batyrev showed that 
the $E$-function of the pair $(\widetilde {M/G},D)$ consisting of 
a resolution $\mu\colon \widetilde {M/G} \rightarrow M/G$ and the divisor 
defined via the discrepancy $D=K_{\widetilde {M/G}}-\mu^*(K_{M/G})$
(with trivial group action) coincides with the orbifold $E$-function. 
The fact that $E$-function of the pair doesn't change under birational
morphisms, as well as an alternative proof of McKay correspondence for
$E$-functions are based on Kontsevich's idea of motivic 
integration 
(cf.\cite{Kontsevich}, \cite{Batyrev}, \cite{denefl}, \cite{looij}).

Another generalization of 
Hirzebruch's $\chi_y$-genus is the (two-variable) elliptic genus,
and this paper grew from an attempt to prove 
the relationship 
between elliptic genera of resolutions of the quotients $M/G$ 
and the elliptic genera associated with the actions of $G$ on $M$.
These two versions of elliptic genus of a global quotient were introduced 
in our previous paper \cite{singellgenus} where the McKay correspondence
was stated as a conjecture.
The proof given below, similarly to Batyrev's approach, requires a  
generalization of elliptic genera considered in \cite{singellgenus}
to the elliptic genus associated with triples consisting of 
a manifold, the group acting on it and the divisor with simple normal 
crossings.

Elliptic genus was extensively studied in recent years (cf.
\cite{Landweber}, \cite{Krichever}, \cite{Hohn}, \cite{Gritsenko}, 
\cite{Totaro} \cite{ellgeninvent}, \cite{MirrorIV} and further references in 
the latter). For an almost complex compact manifold $X$ with
Chern roots $x_i$ (i.e. the total Chern class is $\prod(1+x_i)$)
the elliptic genus can be defined as
\begin{equation}
Ell(X;z,\tau)=\int_X   
\prod_i x_i 
\frac
{\theta (\frac{x_i} {2 \pi \ii}-z,\tau)}
{\theta (\frac{x_i}{2 \pi \ii }, \tau)}
\label{elltheta}
\end{equation} 
where 
$$
\theta(z,\tau)=q^{\frac 18}  (2 \sin \pi z)
\prod_{l=1}^{l=\infty}(1-q^l)
 \prod_{l=1}^{l=\infty}(1-q^l y)(1-q^l y^{-1})
$$ is the 
classical theta function (cf. \cite{Chandra})
where $y=\ee^{2\pi\ii z}$, $q=\ee^{2\pi\ii \tau}$.

Alternatively, the elliptic genus can be written as 
\begin{equation}\label{elleuler}
Ell(X;z,\tau)=\int_X ch(\Ell_{z,\tau})td(X)
\end{equation}
where
$$\Ell_{z,\tau}:=y^{-{\frac {\dim X}2}} \otimes_{n \ge 1} 
\Bigl(
\Lambda_{-yq^{n-1}}
T^*_X \otimes \Lambda_{-y^{-1}q^n} T_{X} \otimes S_{q^n}T^*_X \otimes
S_{q^n} T_X
\Bigr).
$$
Here $T_X$ (resp. $T^*_X$) is complex tangent (resp. cotangent)  
bundle and as usual for a bundle $V$, $\Lambda_t(V)=\sum_i \Lambda^i(V) t^i$
and $S_t(V)=\sum_i Sym^i(V)t^i$ denote 
generating functions for exterior and symmetric 
powers of $V$ (by Riemann Roch this is also the holomorphic Euler 
characteristic of $\Ell_{z,\tau}$). Elliptic genus of a projective
manifold is a holomorphic function
of $(z,\tau) \in {\CC} \times {\cal H}$. Moreover, if $c_1(X)=0$ then 
it is a weak Jacobi form (of weight $0$ and index $\frac{\dim X}2$,
see \cite{ellgeninvent} or earlier references in \cite{MirrorIV}). 

Since 
$y^{-\frac{{\rm dim} X}2}\chi_{-y}(X)=\lim_{q \rightarrow 0} Ell(X;z,\tau)$,
Hirzebruch's $\chi_y$-genus is a specialization of the elliptic genus
(and so are various one variable versions of elliptic genus due to 
Landweber-Stong, Ochanine, Witten and Hirzebruch). On the other 
hand, elliptic genus is a combination of Chern numbers of $X$, 
as is apparent from \eqref{elltheta}, but it cannot be expressed
via Hodge numbers of $X$ (cf. \cite{Hohn}, \cite{ellgeninvent}).
Therefore the information about elliptic genera of resolutions of 
$X/G$ cannot be derived from corresponding information about $E$-function, 
though it can be done for specialization $q \rightarrow 0$ of 
the elliptic genus. Since elliptic genus depends only on the Chern numbers,
it is a cobordism invariant. Totaro \cite{Totaro} 
found a characterization of the elliptic genus \eqref{elltheta}
of $SU$-manifolds from the point of view of cobordisms as 
the universal genus invariant under classical flops. 

A major difference between elliptic genus and $E$-function is that 
the latter is defined 
for quasiprojective varieties. Unfortunately, we don't know if 
a useful definition of elliptic genus can be given for 
arbitrary quasiprojective manifolds. Moreover, while 
$E$-function enjoys strong additivity properties there appears 
to be no analog of them in the case of elliptic genus. Additivity 
allows one to work with $E$-functions not just in the category of 
manifolds but in the category of of arbitrary quasiprojective varieties. 
Nevertheless, in \cite{singellgenus} (extending \cite{ellgeninvent}) 
a definition of elliptic genus for some singular spaces was proposed
as follows.
Let $X$ be a $\QQ$-Gorenstein complex projective variety 
and $\pi\colon Y \rightarrow X$ be a resolution of singularities with 
the simply normal crossing divisor $\cup E_k, k=1,\ldots,r$ as 
its exceptional locus. If the canonical classes of $X$ and $Y$
are related via 
\begin{equation}
K_Y=\pi^*K_X+\sum \alpha_kE_k,
\label{discrepancies}
\end{equation}
then 
\begin{equation}
\widehat{Ell}_Y(X;z,\tau):=
\int_Y 
\Bigl(\prod_l 
\frac{(\frac{y_l}{2\pi\ii})\theta(\frac{y_l}{2\pi\ii}-z)\theta'(0)} 
{\theta(-z)\theta(\frac{y_l}{2\pi\ii})} 
\Bigr)\times 
\Bigl(\prod_k 
\frac{\theta(\frac{e_k}{2\pi\ii}-(\alpha_k+1)z)\theta(-z)} 
{\theta(\frac{e_k}{2\pi\ii}-z)\theta(-(\alpha_k+1)z)} 
\Bigr) 
\label{singellgen}
\end{equation} 
is independent of a resolution $\pi$ (here $e_k$ are the 
cohomology classes of the components $E_k$ of the exceptional divisor
and $y_l$ are the Chern roots of $Y$) and depends only on $X$. 
$\widehat{Ell}_Y(X;z,\tau)$ was called \emph{singular elliptic genus} 
of $X$. When $q \rightarrow 0$, singular elliptic genus specializes
to the singular $\chi_y$-genus calculated from Batyrev's $E$-function.
We refer the reader to \cite{singellgenus} for further discussion of 
this invariant.

On the other hand, for a finite group $G$ of automorphisms of a manifold $X$, 
an \emph{orbifold elliptic genus} was defined in \cite{singellgenus} 
in
terms of the action of $G$ on $X$ as follows.
For a pair of commuting elements $g,h \in G$,
let $X^{g,h}$ be a connected component of the fixed point set of both 
$g$ and $h$. 
Let $TX \vert_{X^{g,h}}=\oplus V_{\lambda}, \lambda(g),
\lambda(h) \in {\QQ} \cap [0,1)$
be decomposition into direct sum, such that $g$ (resp. $h$) 
acts on $V_{\lambda}$ as the 
multiplication by $\ee^{2 \pi \ii \lambda(g)}$
(resp. $\ee^{2 \pi \ii \lambda(h)}$). Then
\begin{equation}
E_{orb}(X,G;z,\tau)=
\label{commuting}
\end{equation}
$$\frac1{\vert G \vert}
\sum_{gh=hg} 
\Bigl(
\prod_{\lambda(g)=\lambda(h)=0}
x_{\lambda}
\Bigr)
\prod_{\lambda} \frac{{\theta(\frac{x_\lambda} {2 \pi \ii }+\lambda(g)-\tau \lambda(h)-z)}
}{
{\theta(\frac{x_\lambda} {2 \pi \ii }+\lambda(g)-\tau \lambda(h))}}
\ee^{2\pi \ii z \lambda(h)z} 
[X^{g,h}].
$$
In \cite{singellgenus} it was conjectured that these two notions of elliptic 
genus coincide. More precisely (cf. Conjecture 5.1, ibid), let $X$ be a 
nonsingular projective variety on which a group $G$ acts effectively 
by biholomorphic transformations. Let $\mu\colon X \rightarrow X/G$ be 
the quotient map and $D=\sum (\nu_i-1) D_i$ be the ramification 
divisor. Let 
$$\Delta_{X/G}:=\sum_j\left(\frac {\nu_j-1}{\nu_j}\right)\mu(D_j).$$
Then 
\begin{equation}\label{maineq}
Ell_{orb}(X,G;z,\tau)=
\left(
\frac 
{2\pi\ii\,\theta(-z,\tau)}
{\theta'(0,\tau)}
\right)^{\dim X}
\widehat{Ell}(X/G,\Delta_{X/G};z,\tau)
\end{equation}
where elliptic genus of the pair $\widehat{Ell}(X/G,\Delta_{X/G};z,\tau)$
is defined by \eqref{singellgen} but with discrepancies $\alpha_k$
obtained from the relation
$$
K_Y=\pi^*(K_{X/G}+\Delta_{X/G})+\sum \alpha_kE_k
$$
rather than the relation \eqref{discrepancies}.

The main goal of this paper is to prove the identity \eqref{maineq} which we
accomplish in Theorem \ref{main}.
One of ingredients of the proof is systematic use of 
the ``hybrid'' orbifold elliptic genus of pairs 
generalizing both singular and orbifold elliptic genera. It is 
defined as follows. Let $(X,E)$ be a Kawamata log-terminal 
pair (cf. \cite{kollar.mori} and Section \ref{section.klt}).
Let $X$ support an action of a finite group $G$ such 
that $(X,E)$ is a $G$-normal pair (cf. \cite{Batyrev} and 
Section \ref{section.maindef}). In addition to notations 
used in the above definition \eqref{commuting} of the orbifold elliptic genus, 
let $\epsilon_k(g), \epsilon_k(h) \in {\QQ} \cap [0,1)$ be defined as
follows. If $E_k$ does not contain $X^{g,h}$ then they are zero
and if $X^{g,h} \subseteq E_k$ then 
$g$ (resp. $h$) acts on ${\cal O}(E_k)$ as the multiplication 
by $\ee^{2 \pi i \epsilon(g)}$ (resp. $\ee^{2 \pi i \epsilon(g)}$).   
Then we define (cf. Definition \ref{maindef}):
\begin{equation}
{\cal Ell}_{orb}(X,E,G;z,\tau):=
\frac 1{\vert G\vert }\sum_{g,h,gh=hg}\sum_{X^{g,h}}[X^{g,h}]
\Bigl(
\prod_{\lambda(g)=\lambda(h)=0} x_{\lambda} 
\Bigr)
\label{intromaindef}
\end{equation}
$$\times\prod_{\lambda} \frac{ \theta(\frac{x_{\lambda}}{2 \pi \ii }+
 \lambda (g)-\tau \lambda(h)-z )} 
{ \theta(\frac{x_{\lambda}}{2 \pi \ii }+
 \lambda (g)-\tau \lambda(h))}  \ee^{2 \pi \ii \lambda(h)z}
$$
$$\times\prod_{k}
\frac
{\theta(\frac {e_k}{2\pi\ii}+\epsilon_k(g)-\epsilon_k(h)\tau-(\delta_k+1)z)}
{\theta(\frac {e_k}{2\pi\ii}+\epsilon_k(g)-\epsilon_k(h)\tau-z)}
{}
\frac{\theta(-z)}{\theta(-(\delta_k+1)z)} \ee^{2\pi\ii\delta_k\epsilon_k(h)z}.
$$

If $G$ is trivial, this expression yields the elliptic genus 
\eqref{singellgen} if $E=\emptyset$ and the 
version of \eqref{singellgen} for pairs as described
earlier for arbitrary $E$. On the other hand, if $G$ is non-trivial but
$E=\emptyset$, then 
one obtains \eqref{commuting}. Moreover $Ell_{orb}(X,E,G)$ for $q
\rightarrow 0$ specializes into Batyrev's 
$E_{orb}(X,E,G;y,1)$  (cf. \cite{Batyrev}).
Thus defined orbifold elliptic genus of pairs is birationally invariant (cf. 
Section \ref{section.maindef}). In fact, we show that 
the contribution of each pair of commuting elements in the 
above definition is invariant under the blowups with normal
crossing nonsingular $G$-invariant centers which allows us to show 
that the contribution of each pair $(g,h)$ is a birational invariant.

Second main ingredient of the proof is the pushforward formula
for the class in \eqref{intromaindef} for toroidal morphisms. Finally,
we use the results of 
\cite{Abramovich.Wang} to show that 
$X \rightarrow X/G$ can be lifted to a toroidal map $\hat Z \rightarrow Z$
so that in the diagram 
$$
\begin{array}{rccc}
&\hat Z&\to&Z\\
 &\downarrow& &\downarrow\\
 &X&\to& X/G
\end{array}
$$
the vertical arrows are resolutions of singularities.

As was already pointed out, singular (resp. orbifold) elliptic
specializes into some known invariants of singular varieties 
(resp. orbifolds). The simplest corollary of our main theorem 
is obtained in the limit 
$q=0$, $y=1$. 
We see that if $X/G$ 
admits a crepant resolution of singularities (i.e. such that 
in \eqref{discrepancies}, one has $\alpha_k=0$ for any $k$) then
the topological euler characteristic of a crepant resolution 
in given by DHVW formula (\ref{eulercommuting}).
While previous proofs of this relation were based 
on motivic integration (cf.  \cite{Batyrev}, \cite{denefl})
the proof presented here uses only birational geometry 
(but depends on \cite{AKMW} and \cite{Abramovich.Wang}).
Moreover, all the results for $E(u,1)$ in quasiprojective case also
get an alternative proof, independent of motivic integration.

Another corollary is the further clarification of a remarkable 
formula due to Dijkgraaf, Moore, Verlinde and Verlinde. 
It was shown in \cite{singellgenus} that
\begin{equation}
\sum_{n\geq 0}p^n Ell_{orb}(X^n/\Sigma_n;z,\tau)=
\prod_{i=1}^\infty \prod_{l,m}
\frac 1
{(1-p^iy^lq^m)^{c(mi,l)}}
\label{dmvv}
\end{equation}
where $\Sigma_n$ is the symmetric group acting on the product 
of $n$ copies of a manifold $X$ such that $Ell(X)=\sum_{m,l}c(m,l)y^lq^m$.
A formula of such type was first proposed in \cite{DMVV}. The main theorem of 
this paper shows that orbifold elliptic genus in \eqref{dmvv} 
can be replaced by singular elliptic genus. While for general $X$ 
it is not clear how to construct a crepant resolution of the symmetric product
(or other kind of 
resolution leading to a calculation of the singular elliptic genus)
in the case $\dim X=2$ it is well known that Hilbert scheme 
$X^{(n)}$ of subschemes of length $n$ in $X$ yields a crepant 
resolution. A corollary of the main theorem is the the following:

{\bf Corollary \ref{hilbertschemes}.}
Let $X$ be a complex projective surface and 
$X^{(n)}$ be its $n$-th Hilbert scheme. Let $\sum_{m,l}c(m,l)y^lq^n$ 
be the elliptic genus of $X$. 
Then 
$$\sum_{n\geq 0}p^n Ell(X^{(n)};z,\tau)=
\prod_{i=1}^\infty \prod_{l,m}
\frac 1
{(1-p^iy^lq^m)^{c(mi,l)}}.
$$

This is a generalization of results due to G\"ottsche on the 
generating series of $\chi_y$-genera of Hilbert schemes 
(cf. \cite{gottsche}) which one obtains for $q=0$.    
In fact in this paper a substantial generalization of \eqref{dmvv}
is proposed. We are able to extend DMVV formula to symmetric powers of
log-terminal varieties and, more generally, to symmetric powers
of Kawamata log-terminal pairs.

The paper is organized as follows. In Section \ref{section.klt} we
recall the concept of Kawamata log-terminal pairs, to the extent necessary
for our purposes. Section \ref{section.maindef} contains our main definition
of the orbifold elliptic genus of a Kawamata log-terminal pair. We 
prove that it is well-defined, for which we use the full force of 
the machinery of \cite{AKMW}.
In Section \ref{section.tormor} we introduce 
toroidal morphisms between pairs that consist of varieties and 
simple normal crossing divisors on them. Our main result is the description
of the pushforward and pullback in the Chow rings 
in terms of the combinatorics of the conical polyhedral complexes.
In the process we use some combinatorial results related to toric varieties,
which are collected in the Appendix \ref{section.app}.
In Section \ref{section.main} we apply these calculations to prove our
main Theorem \ref{main}.
In Section \ref{section.DMVV} we generalize the second quantization
formula of \cite{DMVV} to the case of Kawamata log-terminal pairs.
Various open questions related to our arguments are collected in Section
\ref{open.section}.

The authors would like to thank Dan Abramovich for helpful discussions and
the proof of important Lemma \ref{reduce}. We also thank Nora Ganter whose 
question focused our attention on the problem of defining orbifold elliptic 
genera for pairs.

\section{Kawamata log-terminal pairs}\label{section.klt}

In this section we present the background material for Kawamata log-terminal
pairs which are a standard tool in minimal model program.
Our main reference is \cite{kollar.mori}. 
\begin{proposition}\cite[Definition 2.25, Notation 2.26]{kollar.mori}
Let $(X,D)$ be a pair where $X$ is a normal variety and 
$D=\sum_ia_iD_i$ is a sum of distinct prime divisors on $X$. We allow 
$a_i$ to be arbitrary rational numbers. Assume that $m(K_X+D)$ is 
a Cartier divisor for some $m>0$. Suppose $f\colon Y\to X$ 
is a birational morphism from a normal variety $Y$. Denote by 
$E_i$ the irreducible exceptional divisors and the proper preimages
of the components of $D$.
Then there are naturally defined rational numbers $a(E_i,X,D)$ such that
$$
K_Y=f^*(K_X+D)+\sum_{E_i}a(E_i,X,D)E_i.
$$
Here the equality holds in the sense that a nonzero 
multiple of the difference is a divisor of a rational function.
The number $a(E_i,X,D)$ is called the \emph{discrepancy} of 
$E_i$ with respect to $(X,D)$ and it depends only on $E_i$, but not on $f$. 
By definition $a(D_i,X,D)=-a_i$ and $a(F,X,D)=0$ for any divisor $F\subset X$
different from all $D_i$.
\end{proposition}

\begin{remark}
In the notation of the above proposition,
we will often call the pair $(Y,-\sum_{E_i}a(E_i,X,D)E_i)$ the 
pair on $Y$ that corresponds to $(X,D)$ or the pullback of
$(X,D)$ by $f$. It is easy to see that for any birational
morphism $g\colon Z\to Y$ 
from a normal variety $Z$ the pullback by $g$ of the pullback of $(X,D)$
by $f$ is equal to the pullback of $(X,D)$ by $f\circ g$.
\end{remark}

\begin{definition}
We call a morphism $f\colon Y\to X$ from a nonsingular variety $Y$ to
a normal variety $X$ a resolution of singularities of the pair $(X,D)$
if the exceptional locus of $f$ is a divisor with simple normal crossings,
which is additionally simple normal crossing with the proper
preimage of $D$. Every pair admits a resolution, 
see \cite[Theorem 0.2]{kollar.mori}.
\end{definition}

\begin{definition}
A pair $(X,D)$ is called \emph{Kawamata log-terminal} if there is a resolution
of singularities $f\colon Y\to X$ of $(X,D)$ such that the pullback
$(Y,-\sum_i \alpha_iE_i)$ satisfies $\alpha_i>-1$ for all $i$.
\end{definition}

\begin{remark} 
It is easy to see that our definition of Kawamata log-terminal pair coincides
with \cite[Definition 2.34]{kollar.mori} in view of 
\cite[Corollary 2.31]{kollar.mori}. This corollary also implies that
\emph{any} resolution of singularities of a Kawamata log-terminal 
pair satisfies the condition $\alpha_i>-1$ for all $i$.
\end{remark}

We will also need to describe the behavior of Kawamata log-terminal 
pairs under finite morphisms, in particular under quotient morphisms.
We will use the following result.
\begin{proposition}\cite[Proposition 5.20]{kollar.mori}
Let $g\colon X'\to X$ be a finite morphism between normal varieties.
Let $D'$ and $D$ be $\QQ$-Weil divisors on $X'$ and $X$ respectively such that
$$
K_X'+D'=g^*(K_X+D).
$$
Then $K_X'+D'$ is $\QQ$-Cartier if and only if $K_x+D$ is. 
Moreover, $(X',D')$ is Kawamata log-terminal if and only if $(X,D)$ is.
\end{proposition}

\begin{definition}\label{quotpair}
Let $G$ be a finite group which acts effectively on a normal variety
$X$ and preserves a $\QQ$-Weil divisor $D$. Let $g\colon X\to X/G$ 
be the quotient morphism. Then there is a unique divisor $D/G$ on $X/G$ 
such that 
$$g^*(K_{X/G}+D/G)=K_X+D.$$
The components of $D/G$ are the images of the components of $D$ and 
the images of the ramification divisors of $f$.
We call the pair $(X/G,D/G)$ the \emph{quotient of $(X,D)$ by $G$.}
By the above proposition, the quotient pair is Kawamata log-terminal
iff $(X,D)$ is Kawamata log-terminal.
\end{definition}

We remark that this definition is contained in \cite{Batyrev} in 
the particular case of a smooth variety $X$ and trivial divisor $D$.
It allows us to generalize the definition of the pullback of 
a pair to the case of $G$-equivariant morphisms as follows.

\begin{definition}
Let $g\colon X'\to X$ be a generically finite morphism from a normal
$G$-variety $X'$ to a normal variety $X$ which is birationally 
equivalent to the quotient morphism $f\colon X'\to X'/G$.
We  say that a pair $(X',D')$ is a pullback of a pair $(X,D)$ 
if the pullback of $(X,D)$ to $X'/G$ coincides with the quotient of 
$(X',D')$ by $G$. Just as in the birational case, this pullback 
preserves Kawamata log-terminality.
\end{definition}

\section{Orbifold elliptic genera of pairs}\label{section.maindef}

\begin{definition}\cite{Batyrev}
Let $X$ be a smooth manifold with the action of a finite 
group $G$. Let $E$ be a $G$-invariant divisor on $X$.
The pair $(X,E)$ is called \emph{$G$-normal} if $Supp(E)$ has simple normal
crossings and for every point $x\in X$ the action of the 
isotropy subgroup of $x$ on the set of irreducible components of $Supp(E)$ 
that pass through $x$ is trivial. 
\end{definition}

We will extensively use the theta function $\theta(z,\tau)$
of \cite{Chandra}. By default, the second argument will be $\tau$.
We will suppress it from the notations, unless it is different from $\tau$.
We will implicitly assume the standard properties of $\theta$, namely
its zeroes and transformation properties under the Jacobi group.

\begin{definition}\label{maindef}
Let $(X,E)$ be a Kawamata log-terminal $G$-normal pair
with $E=-\sum_k\delta_kE_k$.
We define \emph{orbifold elliptic class} of the triple $(X,E,G)$ as 
an element of the Chow group $A_*(X)$ by the formula
$$
{\cal Ell}_{orb}(X,E,G;z,\tau):=
\frac 1{\vert G\vert }\sum_{g,h,gh=hg}\sum_{X^{g,h}}[X^{g,h}]
\Bigl(
\prod_{\lambda(g)=\lambda(h)=0} x_{\lambda} 
\Bigr)
$$
$$\times\prod_{\lambda} \frac{ \theta(\frac{x_{\lambda}}{2 \pi \ii }+
 \lambda (g)-\tau \lambda(h)-z )} 
{ \theta(\frac{x_{\lambda}}{2 \pi \ii }+
 \lambda (g)-\tau \lambda(h))}  \ee^{2 \pi \ii \lambda(h)z}
$$
$$\times\prod_{k}
\frac
{\theta(\frac {e_k}{2\pi\ii}+\epsilon_k(g)-\epsilon_k(h)\tau-(\delta_k+1)z)}
{\theta(\frac {e_k}{2\pi\ii}+\epsilon_k(g)-\epsilon_k(h)\tau-z)}
{}
\frac{\theta(-z)}{\theta(-(\delta_k+1)z)} \ee^{2\pi\ii\delta_k\epsilon_k(h)z}.
$$
Here $X^{g,h}$ denotes an irreducible component of the fixed set of the 
commuting elements $g$ and $h$ and $[X^{g,h}]$ denotes the corresponding 
element of $A_*(X)$. The restriction of $TX$ to $X^{g,h}$
splits into linearized bundles according to the ($[0,1)$-valued) characters 
$\lambda$  of $\langle g,h\rangle$, 
see \cite{singellgenus}. In addition, $e_k=c_1(E_k)$ 
and $\epsilon_k$ is the character of ${\cal O}(E_k)$ restricted to $X^{g,h}$ 
if $E_k$ contains $X^{g,h}$ and is zero otherwise.

We define \emph{orbifold elliptic genus} ${Ell}_{orb}(X,E,G)$ of 
$(X,E,G)$ as the degree of the top component of the orbifold elliptic class
${\cal Ell}_{orb}(X,E,G)$.
\end{definition}

\begin{remark}
Notice that in the particular cases of $\vert G\vert =1$ and $E=0$ the above 
definition restricts to that of the singular elliptic genus (up to 
a normalization factor) and orbifold elliptic genus, see \cite{singellgenus}.
However, the notion of orbifold elliptic class appears to be new.
\end{remark}

\begin{remark}
The Kawamata log-terminality assures that we never divide by zero in 
the above formulas.
\end{remark}

Our first goal is to show that orbifold elliptic class is compatible 
with blowups.
\begin{theorem}\label{orbblowup}
Let $(X,E)$ be a Kawamata log-terminal $G$-normal pair
and let $Z$ be a smooth $G$-equivariant locus in
$X$ which is normal crossing to $Supp(E)$.
Let $f:\hat X\to X$ denote the blowup of $X$ along $Z$. We define 
$\hat E$ by $\hat E=-\sum_{k}\delta_k\hat E_k-\delta Exc(f)$
where $\hat E_k$ is the proper transform of $E_k$ and 
$\delta$ is determined from $K_{\hat X}+\hat E = f^*(K_X + E)$.
Then $(\hat X,\hat E)$ is a Kawamata log-terminal $G$-normal pair and 
$$
f_*{\cal Ell}_{orb}(\hat X,\hat E,G;z,\tau)={\cal Ell}_{orb}(X,E, G;z,\tau).
$$
\end{theorem}

\begin{proof}
It is clear that $(\hat X,\hat E)$ is Kawamata log-terminal. Because of the
normal crossing conditions on $Z$ and $Supp(E)$, the divisor $Supp(\hat E)$
has simple normal crossings. The $G$-normality is clearly preserved since
the exceptional divisors do not intersect and any intersection of $\hat E_k$
on $\hat X$ induces an intersection of $E_k$ on $X$.

We will prove the theorem by showing that for every pair $(g,h)$ and every
connected component $X^{g,h}$ the contributions to 
$f_*{\cal Ell}_{orb}(\hat X,\hat E,G;z,\tau)$ 
of connected components $\hat X^{g,h}$
such that $f(\hat X^{g,h})\subseteq X^{g,h}$ equals 
the contribution of $X^{g,h}$ to ${\cal Ell}_{orb}(X,E,G;z,\tau)$.
So from now on $g$, $h$ and $X^{g,h}$ are fixed.

The set of connected components of the fixed point set of $\langle g,h\rangle$ 
that maps inside $X^{g,h}$ is described as follows. Let $Z^{g,h}$ denote
the intersection of $X^{g,h}$ and $Z$. Since $Z$ is $G$-equivariant, the
intersection is a union of some connected components 
of $\langle g,h\rangle$-invariant points of $Z$. Locally at every point 
of the intersection, $Z$ and $X^{g,h}$ intersect normally, since the 
normal spaces to $Z^{g,h}$ inside $Z$ and $X^{g,h}$ have different characters.
For simplicity, we assume that $Z^{g,h}$ is connected, and we will
remark later on the general case.

If $X^{g,h}\neq Z^{g,h}$ then one of the $\hat X^{g,h}$ will be obtained 
as the proper preimage of $X^{g,h}$ under $f$ and will be isomorphic 
to the blowup of $X$ along $Z^{g,h}$. Other components will lie in
the preimage of $Z^{g,h}$ and are described as follows. The restriction
of the normal bundle to $Z$ in $X$ to $Z^{g,h}$ splits into character 
subbundles. For each character $\Lambda$ the projectivization of the 
corresponding bundle over $Z^{g,h}$ is naturally embedded 
into the preimage of $Z^{g,h}$ under $f$ (which is the projectivization
of the whole normal bundle to $Z$ restricted to $Z^{g,h}$).

We first concentrate on the case $X^{g,h}\neq Z^{g,h}$.

Let $N_1$ be the subbundle of the normal bundle to $Z^{g,h}$ in $X$
that is the image of the normal bundle of $Z^{g,h}$ in $Z$.
Let $N_2$ be the subbundle of the normal bundle to $Z^{g,h}$ that is 
the image of the normal bundle of $Z^{g,h}$ in $X^{g,h}$. Finally
let's $N_3$ be the quotient of $NZ^{g,h}$ by the sum of $N_1$ and $N_2$.
The transversality implies that it is also a bundle, i.e. the rank of the 
fibers is constant.

Let us calculate the contribution to $f_*{\cal Ell}_{orb}
(\hat X,\hat E,G;z,\tau)$
that comes from $\hat X^{g,h}_0$, which is the proper preimage of 
$X^{g,h}$, provided $N_2\neq 0$.
As in \cite{singellgenus}, we make a technical assumption that all bundles we 
consider are restrictions of some bundles defined on $X$. The calculation 
follows closely those of \cite{singellgenus}. 
We have
$$
c(T\hat X) = c(f^*TX)(1+\hat z)\Prod_i\frac{(1+f^*m_i-\hat z)}{(1+f^*m_i)}
$$
where $\hat z$ is the first Chern class of the exceptional divisor of $f$
and $\Prod_i(1+m_i)$ is the Chern class of the bundle on $X$ whose 
restriction to $Z$ 
is the normal bundle of $Z$ in $X$. Similarly,
$$
c(T\hat X^{g,h}_0)=c(f^*TX^{g,h})(1+\hat z)\Prod_i
\frac{(1+f^*s_i-\hat z)}{(1+f^*s_i)}
$$
where $\hat z$ and $f$ are restrictions to $X^{g,h}_0$ 
(mild abuse of notation) and $\Prod_i(1+s_i)$ restricts to $c(N_2)$
on $Z^{g,h}$. 

So the Chern class of the normal bundle to $\hat X^{g,h}_0$ is 
$$
c(N\hat X^{g,h}_0) = c(f^*NX^{g,h})\Prod_i 
\frac{(1+f^*t_i-\hat z)}{(1+f^*t_i)}
$$
where $\Prod_i(1+t_i)$ restricts to $c(N_3)$ on $Z^{g,h}$.

We will also need to know how $E_i$ change. For $E_i$ that do not
contain $Z$ we have $\hat E_i=f^*E_i$, and for $E_i$ that contain
$Z$ we have $\hat E_i=f^*E_i-\hat Z$.

As a result, the contribution of $\hat X^{g,h}_0$ to
${\cal Ell}_{orb}(\hat X,\hat E,G;z,\tau)$ is
$$
[X^{g,h}_0] \Prod_i 
\frac {x_i\theta(\frac {x_i}{2\pi\ii}-z)}{\theta(\frac {x_i}{2\pi\ii})}
\prod_{N_2}
\frac{(f^*n_i-\hat z)\theta(\frac{f^*n_i-\hat z}{2\pi\ii}-z)}
{\theta(\frac{f^*n_i-\hat z}{2\pi\ii})}
\frac{\theta(\frac{f^*n_i}{2\pi\ii})}
{(f^*n_i)\theta(\frac{f^*n_i}{2\pi\ii}-z)}
$$
$$\times
\frac{(\frac{\hat z}{2\pi\ii})\theta(\frac{\hat z}{2\pi\ii}-z)\theta'(0)}
{\theta(\frac{\hat z}{2\pi\ii})\theta(-z)}
$$
$$\times
\Prod_{N_1}
\frac{\theta(\frac{f^*n_i}{2\pi\ii}+\lambda_i(g)-\lambda_i(h)\tau-z)}
{\theta(\frac{f^*n_i}{2\pi\ii}+\lambda_i(g)-\lambda_i(h)\tau)}
\ee^{2\pi\ii\lambda_i(h)z}
$$
$$\times
\Prod_{N_3}
\frac{\theta(\frac{f^*n_i-\hat z}{2\pi\ii}+\lambda_i(g)-\lambda_i(h)\tau-z)}
{\theta(\frac{f^*n_i-\hat z}{2\pi\ii}+\lambda_i(g)-\lambda_i(h)\tau)}
\ee^{2\pi\ii\lambda_i(h)z}
$$
$$\times
\Prod_{E_i\supset Z}
\frac
{\theta(\frac{f^*e_i-\hat z}{2\pi\ii}+\epsilon_i(g)-\epsilon_i(h)\tau-(\delta_i+1)z)}
{\theta(\frac{f^*e_i-\hat z}{2\pi\ii}+\epsilon_i(g)-\epsilon_i(h)\tau-z)}
\frac{\theta(-z)}{\theta(-(\delta_i+1)z)}
\ee^{2\pi\ii\delta_i\epsilon_i(h)z}
$$
$$\times
\Prod_{E_i\not\supset Z}
\frac
{\theta(\frac{f^*e_i}{2\pi\ii}+\epsilon_i(g)-\epsilon_i(h)\tau-(\delta_i+1)z)}
{\theta(\frac{f^*e_i}{2\pi\ii}+\epsilon_i(g)-\epsilon_i(h)\tau-z)}
\frac{\theta(-z)}{\theta(-(\delta_i+1)z)}
\ee^{2\pi\ii\delta_i\epsilon_i(h)z}
$$
$$\times
\frac{\theta(\frac{\hat z}{2\pi\ii}-(\delta+1) z)} 
{\theta(\frac{\hat z}{2\pi\ii}-z)} 
\frac{\theta(-z)}{\theta(-(\delta+1) z)}.
$$
In the above formula first two lines account for the tangent bundle
to $\hat X^{g,h}_0$, next two lines account for the normal bundle to it,
and the remaining three lines account for the divisors. We use the notation
$\Prod_{N_i}$ to indicate the product over the Chern roots of the corresponding
bundle.
Notice the normalization factor in the second line.

As in \cite{singellgenus}, we rewrite the above expression as a power 
series $\sum_nR_n\hat z^n$ in $\hat z$. Clearly, $f_*R_0$ is precisely
the contribution of the $X^{g,h}$ to ${\cal Ell}_{orb}(X,E,G;z,\tau)$.
If we denote $r= rk N_2$, we have $f_*\hat z^{r+n}= i_*(s_n(i^*N_2))
(-1)^{n+r-1}$ where $i_*$ is the pushforward from $Z^{g,h}$ to $X^{g,h}$.
We can therefore rewrite the contribution of $f_*R_{>0}$ as
$$
[Z^{g,h}] \sum_{n\geq 0} s_n(i^*N_2)(-1)^{n+r-1}
({\rm Coeff.~at~}\hat z^{r+n})({\rm above~expression}).
$$
Taking into account 
$$\sum_{n\geq 0}s_ni^*N_2(-1)^nt^{-n}=\frac{t^r}{\Prod_{N_2}(t-n_i)},
$$
we can rewrite this as 
$$
(-1)Res_{t=0} [Z^{g,h}]
\frac {\theta(\frac t{2\pi\ii}-(\delta+1)z)\theta'(0)}
{(2\pi\ii)\theta(\frac t{2\pi\ii})\theta(-(\delta+1)z)}
~\Prod_{TZ^{g,h}}\frac{y_i\theta(\frac {y_i}{2\pi\ii}-z)}
{\theta(\frac{y_i}{2\pi\ii})}
$$
$$\times
\Prod_{N_1}
\frac{\theta(\frac{f^*n_i}{2\pi\ii}+\lambda_i(g)-\lambda_i(h)\tau-z)}
{\theta(\frac{f^*n_i}{2\pi\ii}+\lambda_i(g)-\lambda_i(h)\tau)}
\ee^{2\pi\ii\lambda_i(h)z}
~\Prod_{N_2}
\frac{\theta(\frac{f^*n_i-t}{2\pi\ii}-z)}
{\theta(\frac{f^*n_i-t}{2\pi\ii})}
$$
$$\times
\Prod_{N_3}
\frac{\theta(\frac{f^*n_i-t}{2\pi\ii}+\lambda_i(g)-\lambda_i(h)\tau-z)}
{\theta(\frac{f^*n_i-t}{2\pi\ii}+\lambda_i(g)-\lambda_i(h)\tau)}
\ee^{2\pi\ii\lambda_i(h)z}
$$
$$\times
\Prod_{E_i\supset Z}
\frac
{\theta(\frac{f^*e_i-t}{2\pi\ii}+\epsilon_i(g)-\epsilon_i(h)\tau-(\delta_i+1)z)}
{\theta(\frac{f^*e_i-t}{2\pi\ii}+\epsilon_i(g)-\epsilon_i(h)\tau-z)}
\frac{\theta(-z)}{\theta(-(\delta_i+1)z)}
\ee^{2\pi\ii\delta_i\epsilon_i(h)z}
$$
$$\times
\Prod_{E_i\not\supset Z}
\frac
{\theta(\frac{f^*e_i}{2\pi\ii}+\epsilon_i(g)-\epsilon_i(h)\tau-(\delta_i+1)z)}
{\theta(\frac{f^*e_i}{2\pi\ii}+\epsilon_i(g)-\epsilon_i(h)\tau-z)}
\frac{\theta(-z)}{\theta(-(\delta_i+1)z)}
\ee^{2\pi\ii\delta_i\epsilon_i(h)z}.
$$
We will denote the expression above by $F(t)$, to be thought of as 
a meromorphic function on $\CC$ with values in the Chow group 
$A_*(Z^{g,h})$.

Let's now calculate the contributions from other components $\hat X^{g,h}$
that map inside $X^{g,h}$. As we have discussed earlier, these components
correspond to non-trivial characters $\Lambda$ that are present in $N_3$.
We want to find normal and tangent bundle of 
$X^{g,h}_\Lambda\cong\PP N_\Lambda$ inside $\hat X$. The Chern class of
the tangent bundle 
can be described as the restriction from $\hat X$ of 
$$
\Prod_{N_\Lambda} (1+f^*n_i-\hat z)~\Prod_{TZ^{g,h}}(1+f^*y_i),
$$
so the normal bundle has the Chern class which is a restriction 
of
$$
(1+\hat z)\Prod_{N_1}(1+f^*n_i)\Prod_{N_2\oplus N_3/N_\Lambda}
(1+f^*n_i-\hat z).
$$
Therefore, the contribution of $\hat X^{g,h}_{\Lambda}$ to 
${\cal Ell}_{orb}(\hat X,\hat E,G)$ is 
$$
[\PP N_\Lambda]
\frac {\theta'(0)}{2\pi\ii\theta(-z)}
\Prod_{N_\Lambda}
\frac{(f^*n_i-\hat z)\theta(\frac{f^*n_i-\hat z}{2\pi\ii}-z)}
{\theta(\frac{f^*n_i-\hat z}{2\pi\ii})}
\Prod_{TZ^{g,h}}
\frac{f^*y_i\theta(\frac {f^*y_i}{2\pi\ii}-z)}
{\theta(\frac{f^*y_i}{2\pi\ii})}
$$
$$\times
\frac{\theta(\frac{\hat z}{2\pi\ii}+\Lambda(g)-\Lambda(h)\tau-z)}
{\theta(\frac{\hat z}{
2\pi\ii}+\Lambda(g)-\Lambda(h)\tau)}
\ee^{2\pi\ii\Lambda(h)z}
~\Prod_{N_1}
\frac{\theta(\frac{f^*n_i}{2\pi\ii}+\lambda_i(g)-\lambda_i(h)\tau-z)}
{\theta(\frac{f^*n_i}{2\pi\ii}+\lambda_i(g)-\lambda_i(h)\tau)}
\ee^{2\pi\ii\lambda_i(h)z}
$$
$$\times
\Prod_{N_2\oplus N_3/N_\Lambda}
\frac{\theta(\frac{f^*n_i-\hat z}{2\pi\ii}+(\lambda_i-\Lambda)(g)
-(\lambda_i-\Lambda)(h)\tau-z)}
{\theta(\frac{f^*n_i-\hat z}{2\pi\ii}+(\lambda_i-\Lambda)(g)
-(\lambda_i-\Lambda)(h)\tau)}
\ee^{2\pi\ii(\lambda_i-\Lambda)(h)z}
$$
$$\hskip -.1cm\times
\Prod_{E_i\supset Z}
\frac
{\theta(\frac{f^*e_i-\hat z}{2\pi\ii}+(\epsilon_i-\Lambda)(g)
-(\epsilon_i-\Lambda)(h)\tau-(\delta_i+1)z)
\theta(-z)
}
{\theta(\frac{f^*e_i-\hat z}{2\pi\ii}+(\epsilon_i-\Lambda)(g)
-(\epsilon_i-\Lambda)(h)\tau-z)
\theta(-(\delta_i+1)z)
}
\ee^{2\pi\ii\delta_i(\epsilon_i-\Lambda)(h)z}
$$
$$\times
\Prod_{E_i\not\supset Z}
\frac
{\theta(\frac{f^*e_i}{2\pi\ii}+\epsilon_i(g)-\epsilon_i(h)\tau-(\delta_i+1)z)}
{\theta(\frac{f^*e_i}{2\pi\ii}+\epsilon_i(g)-\epsilon_i(h)\tau-z)}
\frac{\theta(-z)}{\theta(-(\delta_i+1)z)}
\ee^{2\pi\ii\delta_i\epsilon_i(h)z}.
$$
$$\times
\frac{\theta(\frac{\hat z}{2\pi\ii}+\Lambda(g)-\Lambda(h)\tau-(\delta+1)z)} 
{\theta(\frac{\hat z}{2\pi\ii}+\Lambda(g)-\Lambda(h)\tau-z)} 
\frac{\theta(-z)}{\theta(-(\delta+1) z)}
\ee^{2\pi\ii \delta\Lambda(h)z}.
$$
Here we used the fact that the line bundle ${\cal O}(\hat Z)$ 
has character $\Lambda$ on $X^{g,h}_\Lambda$. 
We again expand the integrand in terms of powers of $\hat z$ and use
$f_*\hat z^{l-1+n}=s_n(N_\Lambda)(-1)^{l-1+n}$ where $l=rk(N_\Lambda)$,
to rewrite the pushforward to $X$ of the above as 
$$
(-1)Res_{t=0} [Z^{g,h}]
\frac {\theta(\frac t{2\pi\ii}+\Lambda(g)-\Lambda(h)\tau
-(\delta+1)z)\theta'(0)}
{(2\pi\ii)\theta(\frac t{2\pi\ii}+\Lambda(g)-\Lambda(h)\tau)
\theta(-(\delta+1)z)}
\ee^{2\pi\ii(\delta+1)\Lambda(h)z}
$$
$$\times
~\Prod_{TZ^{g,h}}\frac{y_i\theta(\frac {y_i}{2\pi\ii}-z)}
{\theta(\frac{y_i}{2\pi\ii})}
$$
$$\times
\Prod_{N_1}
\frac{\theta(\frac{f^*n_i}{2\pi\ii}+\lambda_i(g)-\lambda_i(h)\tau-z)}
{\theta(\frac{f^*n_i}{2\pi\ii}+\lambda_i(g)-\lambda_i(h)\tau)}
\ee^{2\pi\ii\lambda_i(h)z}
$$
$$\times\Prod_{N_2}
\frac{\theta(\frac{f^*n_i-t}{2\pi\ii}+\Lambda(g)-\Lambda(h)\tau-z)}
{\theta(\frac{f^*n_i-t}{2\pi\ii}+\Lambda(g)-\Lambda(h)\tau)}
\ee^{-2\pi\ii\Lambda(h)z}
$$
$$\times
\Prod_{N_3}
\frac{\theta(\frac{f^*n_i-t}{2\pi\ii}+(\lambda_i-\Lambda)(g)
-(\lambda_i-\Lambda)(h)\tau-z)}
{\theta(\frac{f^*n_i-t}{2\pi\ii}+(\lambda_i-\Lambda)(g)
-(\lambda_i-\Lambda)(h)\tau)}
\ee^{2\pi\ii(\lambda_i-\Lambda)(h)z}
$$
$$\times
\Prod_{E_i\supset Z}
\frac
{\theta(\frac{f^*e_i-t}{2\pi\ii}+(\epsilon_i-\Lambda)(g)
-(\epsilon_i-\Lambda)(h)\tau-(\delta_i+1)z)
\theta(-z)
}
{\theta(\frac{f^*e_i-t}{2\pi\ii}+(\epsilon_i-\Lambda)(g)
-(\epsilon_i-\Lambda)(h)\tau-z)
\theta(-(\delta_i+1)z)
}
\ee^{2\pi\ii\delta_i(\epsilon_i-\Lambda)(h)z}
$$
$$\times
\Prod_{E_i\not\supset Z}
\frac
{\theta(\frac{f^*e_i}{2\pi\ii}+\epsilon_i(g)-\epsilon_i(h)\tau-(\delta_i+1)z)
\theta(-z)
}
{\theta(\frac{f^*e_i}{2\pi\ii}+\epsilon_i(g)-\epsilon_i(h)\tau-z)
\theta(-(\delta_i+1)z)
}
\ee^{2\pi\ii\delta_i\epsilon_i(h)z}
$$
which can be rewritten as 
$$
(-1)Res_{t=\Lambda(g)-\Lambda(h)\tau} F(t)
$$
because the additional exponential factors cancel due to
$\delta=\sum_{E_i\supset Z}\delta_i + rk(N_2)+rk(N_3)-1$.

So in the case $X^{g,h}\neq Z^{g,h}$ all we need is to show that
$$
Res_{t=0}F(t)+\sum_{\Lambda}Res_{t=\Lambda(g)-\Lambda(h)\tau}F(t)=0.
$$
This follows from the observation that $F$ is periodic with respect
to $t\to t+2\pi\ii$ and $t\to t+2\pi\ii\tau$ and has poles at 
$0$ and $\Lambda(g)-\Lambda(h)$ only. Indeed, the periodicity is 
a corollary of the transformation properties of $\theta$  and the  
definition of $\delta$. The statement on poles follows from the 
fact that for every $E_i\supset Z$ the theta function in the denominator
is precisely offset by of the theta functions in the numerator.
Indeed, in view of the normal crossing condition on $Supp(E)$
and $Z$, each $E_k$ gives a quotient bundle of the normal bundle to $Z$
and the sum over all $E_k$ is (locally) a quotient of $N_2\oplus N_3$.
As a result, $e_k$ is a  Chern root of $N_3$ or $N_2$ depending
on whether or not $E_k$ contains $X^{g,h}$.

As in \cite{singellgenus}, 
we remark that we can ignore the assumption that $N_i$ come from bundles
on $X$, because the expression for $F(t)$ makes sense without it and 
deformation to the normal cone can be used in general. We also observe
that in the case when $Z^{g,h}$ has several connected components, the 
above calculation shows that the contributions of the components,
other than $X^{g,h}_0$ to $f_*{\cal Ell}_{orb}(\hat X,\hat E,G;z,\tau)$
cancel the $f_*R_{>0}$ contributions of the connected component $X^{g,h}_0$.
The $f_*R_0$ contribution of $X^{g,h}_0$ is again the contribution of
$X^{g,h}$ to ${\cal Ell}_{orb}(X,E,G;z,\tau)$.

The case $X^{g,h}=Z^{g,h}$ is handled similarly. This time,
the contributions to $Ell(\hat X,\hat E,G;z,\tau)$ equal
$$
-\int_{Z^{g,h}}
\sum_{\Lambda} Res_{t=\Lambda(g)-\Lambda(h)\tau}F(t)
=\int_{X^{g,h}}Res_{t=0}F(t)
$$
which is precisely the 
contribution of $X^{g,h}$ to $Ell(X,E,G;z,\tau)$.
Indeed, since $N_2=0$, and no divisor $E_i$ that contains $Z$ can
have $\epsilon=0$, 
$F(t)$ has a simple pole at $t=0$ and the residue is easy to
calculate.
\end{proof}

We will now use the invariance under blowups to define the orbifold
elliptic genus and orbifold elliptic class for an arbitrary $G$-equivariant
Kawamata log-terminal pair.

\begin{definition}\label{genorb}
Let $(Z,D)$ be an arbitrary $G$-equivariant Kawamata log-terminal pair
with no additional conditions on its singularities. Let 
$\pi:X\to Z$ be a $G$-equivariant 
resolution of singularities of $(Z,D)$, such that 
the corresponding pair $(X,E)$ is $G$-normal. Then the 
orbifold elliptic class of $(Z,D)$ in $A_*(Z)$ is defined as the 
pushforward $\pi_*$ of the orbifold elliptic class of $(X,E)$
and the orbifold elliptic genus of $(Z,D)$ is defined as 
the orbifold elliptic genus of $(X,E)$ or alternatively as the
degree of the orbifold elliptic class.
\end{definition}

Clearly, this definition does not make sense unless we can prove that
it does not depend on the resolution $\pi$.

\begin{theorem}
Definition \ref{genorb} makes sense, that is the pushforwards of the
orbifold elliptic classes do not depend on the resolution of singularities.
\end{theorem}

\begin{proof}
In view of Theorem \ref{orbblowup}, it is enough to show that any 
two $G$-normal resolutions of singularities $(X_-,E_-)$ and $(X_+,E_+)$
of $(Z,D)$ can be connected by a sequence of equivariant 
blowups and blowdowns among $G$-normal resolutions of singularities of 
$(Z,D)$.
This is a $G$-normality strengthening of the equivariant version of
the Weak Factorization Theorem of \cite{AKMW}. The equivariant version
itself assures that such sequence of blowups and blowdowns exists
in the category of simple normal crossing $G$-equivariant divisors $E$.

In order to get $G$-normality, observe that for every simple normal 
crossing $G$-equivariant divisor $E$ on smooth $X$ there is a 
canonical sequence of blowups that makes the preimage $G$-normal.
Namely, this is the toroidal morphism that corresponds to the 
barycentric subdivision of the corresponding polyhedral complex (see 
Section 5.6 of \cite{AKMW}). In the notations of Section 4.3 of \cite{AKMW},
we apply this procedure in the definition of $W_{i\pm}^{res}$. Then
the additional sequences of blowups $r_{i\pm}$ preserve $G$-normality
and the statement is reduced to the case of the toroidal birational
map $\varphi_i^{can}$. The group $G$ acts by interchanging the vertices of 
the polyhedral complexes $\Delta_\pm$ of $W_{i\pm}^{can}$. We apply the
barycentric subdivision blowup to both of them, and then observe that
all intermediate varieties in the toroidal version of weak factorization
have $G$-normal divisors. Indeed, each of them comes from a subdivision
$\Delta$ of $B\Delta_+$ or $B\Delta_-$, where $B$ stands for barycentric
subdivision, and we assume the former with no loss
of generality. If a cone $C$ in $\Delta$ maps to itself by some group 
element $g\in G$, then the same is true for the smallest cone $C_+$
in $B\Delta_+$
that contains its image. However as observed in Section 5.6 of \cite{AKMW},
this implies that $g$ acts trivially on the span of $C^*$, hence on $C$.
This implies $G$-normality, since every fixed point of $g$ comes from
a stratum that corresponds to some cone of $\Delta$.
\end{proof}

\begin{remark}\label{gstrict}
The Weak Factorization Theorem also works in the category of $G$-strict 
divisors, defined by the condition that the translates of every 
irreducible component of $E$ are either equal or disjoint. Indeed, the
above argument works, since $G$-strictness is preserved under 
normal crossing $G$-equivariant blowups with smooth centers and the 
barycentric subdivision assures $G$-strictness, not just $G$-normality.
\end{remark}

\begin{remark}
It is clear from the definition that the orbifold elliptic genus of a 
log-terminal $G$-variety is unchanged under equivariant crepant morphisms.
\end{remark}

\begin{remark}\label{eachgh}
The arguments of this section clearly show that the contribution of
each pair $(g,h)$ of commuting elements of $G$ to the orbifold elliptic
class and genus is well-defined. Indeed, in the proof of 
Theorem \ref{orbblowup}
each pair was considered separately.
\end{remark}

\begin{remark}
Orbifold elliptic genus for the product of triples 
$(X_1,E_1,G_1)$ and $(X_2,E_2,G_2)$ equals the product of elliptic 
genera. The product of the triples is defined as the product of 
the varieties, the sum of the pullbacks of the divisors and the 
direct product of the group actions.
\end{remark}

We observe that our definition of orbifold elliptic genus is compatible
with the definition of the orbifold string $E$-function of 
$E_{orb}(X,E,G)$ of \cite{Batyrev} in the sense that the limit of 
orbifold elliptic genus as $\tau\to \ii\infty$ recovers the orbifold
string function analog of the $\chi_y$-genus. 
For this, we will need the following easy lemma.
\begin{lemma}\label{chiyfree}
Let $X$ be a complete stratified $G$-variety with at most quotient 
singularities such that the action of $G$ is effective and free
and preserves the stratification. Let $X_1$ be any stratum of 
$X$ and let $G_1$ be a subgroup of $G$ that maps $X_1$ to itself.
Then 
$$
\chi_y(X_1/G_1) = \frac 1{|G_1|}\chi_y (X_1).
$$
\end{lemma}

\begin{proof}
We will argue by induction on the dimension of the stratum. In
dimension zero the freeness of the action implies $|G_1|=1$
and $\chi_y(X_1/G_1)=\chi_y(X_1)=1$. For the induction step,
it is enough to assume that $\overline {X_1}=X$ and $X$ is connected.
It is easy to see that the induction assumption allows us to 
consider $X_1$ to be a part of a nonsingular locus of $X$. After 
an equivariant desingularization, we may assume that $X$ is smooth
and $X_1$ is the open stratum. Notice that desingularization
preserves the freeness of the action, which implies
$$
\chi_y(\overline{X/G}) = \frac 1{|G|}\chi_y (\overline {X}).
$$ 
By additivity of $\chi_y$, we can split the above identity according
to the contributions of the strata. Each stratum 
$Y_1$ in $\overline{X/G}$ is a quotient of a stratum $Y$ in $X$. 
If $H$ is a subgroup of $G$ that fixes $Y$, then there are $|G:H|$ disjoint 
strata of $X$ that map to $Y_1$. By induction assumption, $\chi_y(Y_1)
=\frac 1{|H|}\chi_y(Y)=\frac 1{|G|}\sum_{\{gY\}}\chi_y(gY)$ where the
sum is taken over the cosets of $H$. Consequently, the terms corresponding
to smaller dimensional strata cancel, which finishes the proof of the lemma.
We remark that the statement generally fails for 
free actions on the noncomplete varieties. It is crucial that the 
action stays free on the completion of the stratum.
\end{proof}

\begin{proposition}\label{compbat}
Let $E_{orb}(X,E,G;u,v)$ be defined as in \cite{Batyrev}. Then
$$
\lim_{\tau\to\ii\infty}Ell_{orb}(X,E,G;z,\tau)=y^{-\frac {\dim X}2}
E_{orb}(X,E,G;y,1)$$
where $y=\ee^{2\pi\ii z}$.
\end{proposition}

\begin{proof}
From the product formula for $\theta$, see \cite{Chandra}, we have
$$\lim_{\tau \to \ii\infty}  
\frac
{\theta(u-\beta,\tau)}
{\theta(u,\tau)}
=
\frac
{(1-\ee^{-2 \pi \ii (u-\beta)})} 
{(1-\ee^{-2 \pi \ii u})} 
\ee^{-\pi \ii\beta}$$
and
$$
\lim_{\tau \to \ii\infty}  
\frac
{\theta(u-\alpha\tau-\beta,\tau)}
{\theta(u-\alpha\tau,\tau)}
=\ee^{-\pi \ii\beta}$$
for $0<\alpha<1$.
Hence, by taking the limit in Definition \ref{maindef},
$$\lim_{\tau \to \ii\infty} Ell_{orb}(X,E,G;z,\tau)=
\frac 1 {\vert G \vert}\sum_{g,h, gh=hg} \sum_{X^{g,h}} 
\int_{X^{g,h}}\prod_{\lambda(g)=\lambda(h)=0} x_{\lambda}
 \frac
{(1-\ee^{-x_{\lambda}+2 \pi \ii z})}
{(1-\ee^{-x_{\lambda}})}
$$
$$\times\,
\ee^{-\pi\ii (\dim X) z}
\ee^{2 \pi \ii (\sum_\lambda \lambda(h)) z }
\prod_{\lambda(h)=0,\lambda(g)\ne 0}
\frac
{(1-\ee^{-x_\lambda-2 \pi \ii \lambda(g)+2 \pi \ii z})}
{(1-\ee^{-x_\lambda-2 \pi \ii \lambda(g)})}
$$ 
$$\times\,
\ee^{2 \pi \ii \sum_k\delta_k \epsilon_k(h)z}
\prod_{k,\epsilon_k(h)=0}
\frac
{(1-\ee^{-e_k-2 \pi \ii \epsilon_k(g)+2 \pi \ii (\delta_k+1)z})} 
{(1-\ee^{-e_k-2 \pi \ii \epsilon_k(g)+2 \pi \ii z})}
\frac
{(1-\ee^{2 \pi \ii z})} 
{(1-\ee^{2 \pi \ii (\delta_k+1)z})}
$$$$\times
\prod_{k,\epsilon_k(h)\neq 0}
\frac
{(1-\ee^{2 \pi \ii z})} 
{(1-\ee^{2 \pi \ii (\delta_k+1)z})}
$$
$$ =
\frac {y^{-\frac{\dim X}2}}{ \vert G \vert } 
\sum_{g,h, gh=hg} \sum_{X^{g,h}} 
\int_{X^{g,h}}td(X^{g,h}) 
y^{wt(h,X^h,E)}
$$ 
$$\times
\frac
{ch (\Lambda_{-y} \Omega^1_{X^{h}} \vert_{X^{g,h}}(g))}
{ch (\Lambda_{-1}N^*_{X^{g,h} \subseteq X^h}(g))}
\,
\prod_{E_k\not\supseteq X^h}
\frac
{(1-\ee^{-e_k-2 \pi \ii \epsilon_k(g)}y^{\delta_k+1})} 
{(1-\ee^{-e_k-2 \pi \ii \epsilon_k(g)}y)}
\prod_{E_k}
\frac{(1-y)}{(1-y^{\delta_k+1})}. 
$$
Here $wt(h,X^h,E)$ is the same weight as defined in \cite{Batyrev},
for the irreducible component $X^h$ of the fixed point set of 
$h$ that contains $X^{g,h}$.  We have also used
$$
ch (\Lambda_{-y}\Omega^1_{X^h} \vert_{X^{g,h}}(g))
=\prod_{\lambda(h)=0} 
(1-\ee^{-x_{\lambda}-2 \pi \ii \lambda(g)+2 \pi \ii z})
$$ 
and 
$$ch(\Lambda_{-1} N^*_{X^{g,h} \subseteq X^h}(g))=
\prod_{\lambda(h)=0, \lambda(g) \ne 0}
 ({1-\ee^{-x_\lambda-2 \pi \ii \lambda(g)}}).$$
We use a trick to rearrange the product over $E_k\not\supseteq X^h$ as follows.
$$ \lim_{\tau \to \ii\infty} Ell_{orb}(X,E,G;z,\tau)
$$
$$=
\frac {y^{-\frac{\dim X}2}}{ \vert G \vert } 
\sum_{g,h, gh=hg} \sum_{X^{g,h}} 
\int_{X^{g,h}}td(X^{g,h}) 
y^{wt(h,X^h,E)} 
\frac
{ch (\Lambda_{-y} \Omega^1_{X^h} \vert_{X^{g,h}}
(g))} 
{ch (\Lambda_{-1}N^*_{X^{g,h} \subseteq X^h})}$$
$$
\times\prod_{E_k\not\supseteq X^h}
\Big(1+
\frac
{(y-y^{\delta_k+1})(1-\ee^{-e_k-2 \pi \ii \epsilon_k(g)})}
{(y^{\delta_k+1}-1)(1-y\ee^{-e_k-2 \pi \ii \epsilon_k(g)})}
\Big)
\prod_{E_k\supseteq X^h}
\frac
{(y-1)} 
{(y^{\delta_k+1}-1)}
$$
$$
=\frac {y^{-\frac{\dim X}2}}{ \vert G \vert } 
\sum_{g,h, gh=hg} \sum_{X^{g,h}} 
\int_{X^{g,h}}td(X^{g,h}) y^{wt(h,X^h,E)} 
\frac
{ch (\Lambda_{-y} \Omega^1_{X^{h}} \vert_{X^{g,h}}(g))} 
{ch (\Lambda_{-1}N^*_{X^{g,h} \subseteq X^h})}
$$
$$\times\,
\sum_{J \subseteq I(X^h)} 
\prod_{k \in J}
\frac{(y-y^{\delta_k+1})(1-\ee^{-e_k-2 \pi \ii \epsilon_k(g)})}
{(y^{\delta_k+1}-1)(1-y\ee^{-e_k-2 \pi \ii \epsilon_k(g)})}
\prod_{E_k\supseteq X^h}
\frac
{(y-1)} 
{(y^{\delta_k+1}-1)}
$$
$$
=\frac {y^{-\frac{\dim X}2}}{ \vert G \vert } 
\sum_{g,h, gh=hg} \sum_{X^{g,h}} 
\sum_{J \subseteq I(X^h)} 
\int_{X^{g,h}\cap E_J}td(X^{g,h}\cap E_J) y^{wt(h,X^h,E)} 
$$
$$\times\,
\frac
{ch (\Lambda_{-y} \Omega^1_{X^{h} \cap E_{J}}) \vert_{X^{g,h} \cap E_J}(g)}
{ch(\Lambda_{-1}N^*_{X^{g,h} \cap E_J \subseteq X^h \cap E_J})(g)}
\prod_{k \in J}
\frac{(y-y^{\delta_k+1})}
{(y^{\delta_k+1}-1)}
\prod_{E_k\supseteq X^h}
\frac
{(y-1)} 
{(y^{\delta_k+1}-1)}
.
$$
Here the set $I(X^h)$ is defined as the set of all $k$, such 
that $E_k\not\supseteq X^h$ and $E_k\cap X^h\neq \emptyset$,
which in particular implies that $E_k$ is mapped to itself by $h$
due to $G$-normality.
We have also used the identities
$$td(X^{g,h} \cap E_J) 
=td(X^{g,h})
{\prod_{k \in J, X^{g,h} \not\subseteq E_k}}
\frac {(1-\ee^{-{e_k}})}{e_k} 
$$ 
$$ch (\Lambda_{-y} \Omega^1_{X^{h} \cap E_{J}}) \vert_{X^{g,h} \cap E_J}(g)
=
ch (\Lambda_{-y} \Omega^1_{X^{h}}) \vert_{X^{g,h}}(g)
\prod_{k \in J}
(1-y\ee^{-e_k-2 \pi \ii \epsilon_k(g)})^{-1}
$$
$$ 
ch(\Lambda_{-1}N^*_{X^{g,h} \cap E_J \subseteq X^h \cap E_J})(g)
=
ch (\Lambda_{-1}N^*_{X^{g,h} \subseteq X^h})(g)
\prod_{k\in J,X^{g,h}\subseteq E_k}
\hskip -.1cm
(1-\ee^{-e_k-2 \pi \ii \epsilon_k(g)}).
$$
Changing the order of summation, one obtains
$$ \lim_{\tau \to \ii\infty} Ell_{orb}(X,E,G;z,\tau)
$$
$$=\frac {y^{-\frac{\dim X}2}}{ \vert G \vert } 
\sum_{h\in G}\sum_{X^h}y^{wt(h,X^h,E)} 
\sum_{J \subseteq I(X^h)} 
\sum_{g \in C(h,X^h,J)}  
\int_{X^{g,h} \cap E_J}td(X^{g,h} \cap E_J) 
$$
$$\times\,
\frac
{ch (\Lambda_{-y} \Omega^1_{X^{h} \cap E_J})\vert_{X^{g,h} \cap E_J}
(g)}
{ch (\Lambda_{-1}N^*_{X^{g,h} \cap E_J\subseteq X^h\cap E_J})(g)}
\prod_{k \in J}(\frac{y-1} 
{y^{\delta_k+1}-1}-1)
\prod_{E_k\supseteq X^h}
\frac
{(y-1)} 
{(y^{\delta_k+1}-1)}
$$ 
where the group $C(h,X^h,J)$ is defined as the subgroup of the centralizer 
of $h$ that consists of group elements that map $X^h$ to itself and preserve 
all elements of $J$. 
By the equivariant Riemann Roch theorem the above expression equals 
$$
\frac {y^{-\frac{\dim X}2}}{ \vert G \vert } 
\sum_{h,X^h,J\subseteq I(X^h)}
y^{wt(h,X^h,E)}  {\vert C(h,X^h,J) \vert} 
\chi_y (X^h \cap E_J/C(h,X^h,J)) 
$$
$$\times
\prod_{k\in J}(\frac{y-1}
{y^{\delta_k+1}-1}-1)
\prod_{E_k\supseteq X^h}
\frac
{(y-1)} 
{(y^{\delta_k+1}-1)}
$$
$$
=\frac {y^{-\frac{\dim X}2}}{ \vert G \vert } 
\sum_{h,X^h,J\subseteq I(X^h)}
y^{wt(h,X^h,E)}{ \vert C(h,X^h,J) \vert} 
\chi_y (X^h \cap E_J/C(h,X^h,J)) 
$$
$$\times
\sum_{J_1\subseteq J}(-1)^{\vert J\vert-\vert J_1\vert}\prod_{k\in J_1}
\frac{(y-1)}
{(y^{\delta_k+1}-1)}
\prod_{E_k\supseteq X^h}
\frac
{(y-1)} 
{(y^{\delta_k+1}-1)}
.$$
We observe that we can replace the group $C(h,X^h,J)$ by a possibly
bigger group $\hat C(h,X^h,J)$ characterized by the condition of
fixing $h$ and $X^h$ and fixing $J$ \emph{as a set}. Indeed, the 
$G$-normality of $E$ implies that the action of 
$ \hat C(h,X^h,J)/C(h,X^h,J)$ on $X^h\cap E_J/C(h,X^h,J)$ is free
and we can rewrite the above as 
$$
\frac {y^{-\frac{\dim X}2}}{ \vert G \vert } 
\sum_{h,X^h,J\subseteq I(X^h)}
y^{wt(h,X^h,E)}{ \vert \hat C(h,X^h,J) \vert} 
\chi_y (X^h \cap E_J/\hat C(h,X^h,J)) 
$$
$$\times
\sum_{J_1\subseteq J}(-1)^{\vert J\vert-\vert J_1\vert}\prod_{k\in J_1}
\frac{(y-1)}
{(y^{\delta_k+1}-1)}
\prod_{E_k\supseteq X^h}
\frac
{(y-1)} 
{(y^{\delta_k+1}-1)}
.$$
The variety $X^h$ 
is stratified by intersections with various $E_J^\circ$
which induces a stratification on $X^h\cap E_J/\hat C(h,X^h,J)$. 
Every $J_2\supseteq J$ gives a stratum $X^h\cap E^\circ_{J_2}$
on $X^h\cap E_J$, but different such strata may map to the same stratum
in $X^h\cap E_J/\hat C(h,X^h,J)$. In fact, the strata for 
all possible sets of $J_2$ from
the same orbit of $\hat C(h,X^h,J)$-action on the set of $J_2$ that contain $J$
will map to the same stratum $S$ on $X^h\cap E_J/\hat C(h,X^h,J)$. This stratum
$S$ will be isomorphic to $X^h\cap E^\circ_{J_2}/\hat C(h,X^h,J\subseteq J_2)$ 
where $\hat C(h,X^h,J\subseteq J_2)$ 
is the subgroup of $G$ that fixes $h$ and $X^h$ and fixes $J$ and $J_2$ 
as sets. By Lemma \ref{chiyfree}, we get 
$$
\chi_y(S)=
\frac{|\hat C(h,X^h,J_2)|}
{|\hat C(h,X^h,J\subseteq J_2)|}
\chi_y(X^h \cap E^\circ_{J_2}/\hat C(h,X^h,J_2)).
$$
Indeed, both groups $\hat C(h,X^h,J\subseteq J_2)$ and $\hat C(h,X^h,J_2)$
act freely on the variety $X^h \cap E_{J_2}/C(h,X^h,J_2)$ and preserve 
the stratification which allows one to compare the $\chi_y$-genera of
the quotients. Using additivity
property of $\chi_y$-genus we now get
$$
\chi_y(X^h \cap E_J/\hat C(h,X^h,J))
=
\sum_{J_2\supseteq J}
\frac
{\vert \hat C(h,X^h,J_2) \vert}
{\vert \hat C(h,X^h,J) \vert}
\chi_y(X^h \cap E^\circ_{J_2}/\hat C(h,X^h,J_2))
$$
with the rational coefficients included to account for the fact that 
the same stratum on the quotient may come from different strata on
$X^h\cap E_J$. We notice that $\sum_{J,J_2\supseteq J\supseteq J_1}
(-1)^{\vert J\vert-\vert J_1\vert}$ equals $1$ for $J_1=J_2$ and equals
zero otherwise, to get
$$
\lim_{\tau \to \ii\infty} Ell_{orb}(X,E,G;z,\tau)
=\frac {y^{-\frac{\dim X}2}}{ \vert G \vert } 
\sum_{h,X^h,J\subseteq I(X^h)}
y^{wt(h,X^h,E)}
$$
$$\times
{ \vert \hat C(h,X^h,J) \vert} 
\chi_y (X^h \cap E^\circ_J/\hat C(h,X^h,J)) 
\prod_{k\in J}
\frac{(y-1)}
{(y^{\delta_k+1}-1)}\prod_{E_k\supseteq X^h}
\frac
{(y-1)} 
{(y^{\delta_k+1}-1)}
$$
$$
= {y^{-\frac{\dim X}2}}
\sum_{\{h\},\{X^h\},\{J\}}
y^{wt(h,X^h,E)}
\chi_y (X^h \cap E^\circ_J/\hat C(h,X^h,J)) 
$$$$\times
\prod_{k\in J}
\frac{(y-1)}
{(y^{\delta_k+1}-1)}
\prod_{E_k\supseteq X^h}
\frac
{(y-1)} 
{(y^{\delta_k+1}-1)}
.$$
Here we are summing over representatives $h$ of conjugacy classes of
$G$, then over representatives $X^h$ of the orbits of the action
of $C(h)$ on the components of the fixed point set of $h$ and finally
over the orbits of the action of $C(h,X^h,\emptyset)$ on the subsets of 
 $I(X^h)$. 
This can be compared with Definitions 6.1 and 6.3 of \cite{Batyrev}.
Our sum over the subsets of the set of components fixed by $h$ that 
contain the set of components $E_k$ that contain $X^h$ coincides with 
the set from the definition of \cite{Batyrev} up to trivial contributions.
Indeed, in Definition 6.1 of \cite{Batyrev} $W_J$ is empty unless
$J$ consists of the elements that correspond to divisors that intersect $W$
and moreover contains all elements that correspond to the divisors that
contain $W$.

However, it appears that we are summing over the orbits $\{J\}$ whereas 
Definition 6.3 of \cite{Batyrev} contains the sum over all $J$. The 
extra factor is equal to the length of the orbit of $J$ under the 
action of $C(h,X^h,\emptyset)$. This appears to be a typo in
\cite{Batyrev}, which can be easily seen for a fixed point free action
of $G$.
\end{proof}

\begin{remark} 
Clearly, the comparison between orbifold elliptic genus and orbifold
$E$-function follows from Theorem \ref{main} and the main result of
\cite{Batyrev}. However, it would be strange to rely on such a roundabout
way of proving it.
\end{remark}

\begin{proposition}
Let $X$ be a smooth $G$-variety and let $E$ be a $G$-normal divisor on it
such that $(X,E)$ is Kawamata log-terminal. Let $m(K_X+E)$ be a trivial
Cartier divisor for some integer $m$. 
Denote by $n$ the order of the image of the homomorphism 
$G \rightarrow Aut H^0(X,m(K_X+E))$ where the 
homomorphism can be defined due to $G$-invariance of $E$.
Then 
$Ell_{orb}(X,E,G)$ is a weak Jacobi form of weight $0$ and index $\dim X/2$
with respect to subgroup of the Jacobi group $\Gamma^J$
generated by transformations
$$(z,\tau) \to (z+mn,\tau),
~(z,\tau) \to (z+mn\tau, \tau),
~(z,\tau)\to (z,\tau +1), 
~(z,\tau)\to (\frac z\tau,-\frac {1}\tau).$$
\end{proposition}

\begin{proof}
As in the proof of Theorem 4.3 in \cite{singellgenus}, we introduce
 $$\Phi(g,h,\kappa,z,\tau,x):=
\frac{\theta(\frac {x} {2 \pi \ii }+\kappa(g)-\tau \kappa(h)-z)}
{\theta(\frac x {2 \pi \ii }+\kappa(g)-\tau \kappa(h))}
\ee^{2\pi \ii z \kappa(h)}$$ 
where $\kappa$ is a character of the subgroup of $G$ generated by
$g$ and $h$ considered acting on a line bundle with the first Chern class 
$x$.
Then the contribution of a connected component $X^{g,h}$ 
in Definition \ref{maindef} 
is 
$$
\Bigl(
\prod_{\lambda(g)=\lambda(h)=0} x_{\lambda} 
\Bigr)
\times\prod_{\lambda} 
\Phi(g,h,\lambda,z,\tau,x_{\lambda})$$ 
$$ \prod_k 
\frac{\Phi(g,h,\epsilon_k,(1+\delta_k)z,\tau,e_k)}  
{\Phi(g,h,\epsilon_k,z,\tau,e_k)}
\frac
{\theta(-z)} 
{\theta(-(\delta_k+1)z)}[X^{g,h}].
$$
The proposition follows from transformation
properties of $\Phi(g,h,\kappa,z,\tau,x)$ proven 
in Theorem 4.3 of \cite{singellgenus}. Note that these properties yield 
that transformation $(z,\tau)
\rightarrow (z+mn\tau)$ transforms $Ell_{orb}(X,E,G)$ 
as a Jacobi form provided: $\sum_{\lambda} x_{\lambda}+
\sum \delta_kE_k=0$ and for any $g \in G$ one has 
$mn(\sum_{\lambda} \lambda(g)+
\sum_k \delta_k \epsilon_k (g)) \in {\bf Z}$.
Those are the assumptions of the proposition. 
Jacobi property for the transformation $(z,\tau) \rightarrow (z+1,\tau)$ 
also uses $mn(\sum_{\lambda} \lambda(g)+
\sum_k \delta_k \epsilon_k (g)) \in {\bf Z}$. The other two generators
of $\Gamma^J$ mentioned above transform the contribution of the pair $(g,h)$ 
into the contribution for pairs $(gh^{-1},h)$ and $(h,g^{-1})$ respectively,
multiplied by corresponding Jacobi factor.

We remark that the result of this proposition also follows from the main
Theorem \ref{main} of this paper and \cite[Proposition 3.8]{singellgenus}.
\end{proof}

\section{Toroidal morphisms of nonsingular pairs}\label{section.tormor}
The goal of this section is to derive pullback and pushforward formulas
for functions of divisor classes for certain maps of varieties with 
simple normal crossing divisors on them. 

Let $Z$ be a smooth algebraic variety, together with an open set $U_Z$
whose complement is a simple
normal crossing divisor $D=\sum_{i\in I_Z} D_i$
where $D_i$ are the irreducible components of $D$.
To every subset $I\subseteq I_Z$ and every connected component $Z_{I;j}$ of 
$Z_I=\cap_{i\in I}D_i$ we associate a cone $C_{I;j}$ in the lattice 
$N_{I;j}\cong \ZZ^{\vert I\vert }$. 
We denote the standard basis of $N_{I;j}$ by $\{e_{k;j}\},k\in I$.
The cone $C_{I;j}$ is defined as $\oplus_{k\in I}\RR_{\geq 0}e_{k;j}$. 
For any cone $C$ its relative interior will be denoted by $C^\circ$.

If $I_1\subseteq I_2$ and a connected component $Z_{I_1;j_1}$ contains 
a connected component $Z_{I_2;j_2}$ then we define a 
\emph{face inclusion map} from $N_{I_1;j_1}$ to $N_{I_2;j_2}$
by mapping $e_{k;j_1}$ to $e_{k;j_2}$ for every $k\in I_1$. The image of 
the cone $C_{I_1;j_1}$ under this map is a face of $C_{I_2;j_2}$, which 
explains the terminology. In agreement with the terminology 
of \cite{KKMS} we define the 
\emph{conical polyhedral complex} $\Sigma_Z$  of $(Z,D)$ as the union 
of all cones $C_{I;j}$ glued according to 
the face inclusion maps. We will often refer to it as conical complex.
This is the same as conical polyhedral complex 
with an integral structure for the smooth toroidal embedding without 
self-intersection in the terminology of \cite{KKMS}. We also observe
that closed subvarieties $Z_{I;j}$ induce a stratification on $Z$.
The corresponding locally closed strata will be denoted by $Z^\circ_{I;j}$.

We define piecewise linear (resp. polynomial) functions on $\Sigma_Z$ 
as collections of linear (resp. polynomial) functions on each 
$C_{I;j}\in \Sigma_Z$ which are compatible with all face inclusions. 
We will analogously talk about formal power series on the conical complex 
by considering the completion of the space of the polynomial functions by 
the degree filtration, i.e. the space of collections of formal power series 
on the vector space $N_{I;j}\otimes_\ZZ \CC$ for each $Z_{I;j}$ that are
compatible with the face inclusions. 
There is a natural ring structure on the space of formal power series, 
which we will denote by $\CC[[\Sigma_Z]]$. 

Another natural ring to consider is the partial semigroup ring defined 
by the conical complex $\Sigma_Z$. It is a vector space whose
basis elements $[v]$ are in one-to-one correspondence with lattice points 
$v$ of $\Sigma_Z$. For every pair of points $v_1,v_2\in \Sigma_Z$,
the product $[v_1][v_2]$ is defined as follows
$$
[v_1][v_2] = \tilde{\sum_{\substack{C\in\Sigma_Z\\v_1,v_2\in C}}}[v_1+v_2].
$$
where $\tilde{\sum}$ means that the same point of $\Sigma_Z$ that appears
from different cones is counted only once. Alternatively, it is 
enough to consider the cones $C\ni v_1,v_2$ that do not contain 
any smaller such cone.
In particular, the product is zero if there are no cones $C$ that contain 
both $v_1$ and $v_2$. This ring will be denoted by $\CC[\Sigma_Z]$. 
It can also be thought of as a subring of the direct sum of the semigroup
rings $\oplus_{I;j}\CC[C_{I;j}]$ that consists of collections that are 
compatible with the face inclusions. The identification is via mapping
$[v]$ to the collection of $[v]$ for $C\ni v$ and $0$ otherwise.

It will be crucial to our calculations to construct a natural isomorphism
between the ring $\CC[\Sigma_Z]$ and the subring of $\CC[[\Sigma_Z]]$ 
that consists of piecewise polynomial functions. Namely, for every 
cone $C_{I;j}$ we denote by $x_{k;j}$ the linear functions on 
$N_{I;j}$  such that $x_{k;j}(e_{l;j})=\delta_k^l$ where $\delta$ is 
the Kronecker symbol. The element $[v]=[\sum_{k\in I}a_ke_{k;j}]$ 
of $\CC[C_{I;j}]$ is mapped to the polynomial $\Prod_{k\in I}(x_{k;j})^{a_k}$.
If a collection of elements of $\CC[C_{I;j}]$ is compatible with face 
restrictions, then so is the collection of the corresponding 
polynomial functions. Indeed for any face inclusion between $C_{I_1;j_1}$
and $C_{I_2;j_2}$ the linear functions $x_{k;j_2}$ restrict to 
$x_{k;j_1}$ if $k\in I_1$ and to $0$ otherwise. It is straightforward to
see that this identification is compatible with the product structure.
The inverse map from piecewise polynomial functions on $\Sigma_Z$ to
$\CC[\Sigma_Z]$ is easy to construct as well. In what follows we will
frequently pass from one description of $\CC[\Sigma]$ to the other.

\begin{definition}
We define a map $\rho\colon\CC[[\Sigma_Z]]\to A^*(Z)$ 
as follows. For every lattice point $v\in C_{I;j}^\circ$
given by 
$$
v=\sum_{i\in I}k_ie_{i;j}, k_i\geq 1
$$
we define by $f$ the corresponding piecewise polynomial function on
$\Sigma_Z$ and set
$$
\rho(f)=Z_{I;j}\cap(\cap_{i\in I}(D_i)^{k_i-1}).
$$
We extend the definition of $\rho$ to arbitrary piecewise polynomial
functions by linearity. We extend it to arbitrary piecewise formal power 
series by noticing that that only $v$ with $\sum_i k_i\leq \dim Z$ contribute
nontrivially.
\end{definition}

\begin{proposition}\label{rhoishomo}
The map $\rho$ defined above is a ring homomorphism.
\end{proposition}

\begin{proof}
It is enough to calculate the image of the product of two 
monomial functions $f_1$ and $f_2$ that correspond to points
$v_1$ and $v_2$ in the conical complex. If there is no cone 
$C_{I;j}\in\Sigma_Z$ that contains both $v_1$ and $v_2$, then 
$f_1f_2=0$. On the other hand, in this case the components 
$Z_{I_1;j_1}$ and $Z_{I_2;j_2}$ do not intersect, so
$\rho(f_1)\rho(f_2)=0$.

In general, the product $f_1f_2$ will correspond to 
$$
\sum_{C_{I;j}\supseteq C_{I_1;j_1},C_{I;j}\supseteq C_{I_2;j_2}}[v_{C_{I;j}}]
$$
where $I=I_1\cup I_2$ and
$$
v_{C_{I;j}}=\sum_{i\in I_1} k_{i,1}e_{k;j} +\sum_{i\in I_2} k_{i,2}e_{k;j}.
$$
The cones $C_{I;j}$  are in one-to-one correspondence with the connected 
components of the intersection $Z_{I_1;j_1}\cap Z_{I_2;j_2}$. The image
of each $f_{C_{I;j}}$ under $\rho$ is  
$$
\rho(f_{C_{I;j}})=Z_{I;j}\cap(\cap_{i\in I_1\cup I_2}D_i^{k_{i,1}+k_{i,2}-1})
$$
where $k_{i,1}$ is defined to be zero for $i\not\in I_1$ and similarly
for $k_{i,2}$. On the other hand, the excess intersection formula
\cite{Fulton} gives 
$$
Z_{I_1;j_1}\cap Z_{I_2;j_2}=
\sum_{j} Z_{I;j}\cap (\cap_{i\in I_1\cap I_2}D_i)
$$
in $A^*(Z)$. Then it is easy to see that $\rho(f_1f_2)=\rho(f_1)\rho(f_2)$.
\end{proof}

\begin{remark} 
If all $Z_I=\cap_{i\in I}D_i$ are either empty or connected, then the 
above discussion simplifies greatly. Then the image of $\rho$ is 
precisely the subring of $A^*(Z)$ generated by the classes of $D_i$.
The difficulty was to somehow localize a polynomial in $D_i$ to the correct 
connected component.
\end{remark}

\begin{remark}
The relation between lattice points of $\Sigma_Z$ and piecewise polynomial
functions on $\Sigma_Z$ becomes important for what follows. While the 
former are easier to describe, the latter behave better under pullbacks.
\end{remark}

We now define a certain class of morphisms between two varieties $\hat Z$ 
and $Z$ with the normal crossing divisors $\hat D$ and $D$ respectively.
This is a particular case of the general definition of \cite{Abramovich.Karu}.
\begin{definition}\label{tormor}
We call a proper generically finite 
morphism $\mu\colon \hat Z\to Z$ \emph{toroidal} 
if the following conditions hold.

$\bullet$ $\hat D=\mu^{-1}D$ and the morphism $\mu$ is finite 
and nonramified outside of $\hat D$.

$\bullet$ The image of the closure of  any stratum of $\hat Z$ is the 
closure of a stratum in $Z$.

$\bullet$ For every points  $\hat z\in \hat Z$ and $z\in Z$ 
such that $\mu(\hat z)=z$ and every system of local analytic coordinates
at $z$ such that components of $D$ that pass through $z$ are coordinate
hyperplanes, there exists a system of local analytic coordinates at $\hat z$ 
such that the map $\mu$ is given by monomials.

\end{definition}

\begin{figure}[tbh]
\begin{center}
\includegraphics[scale = 0.4]{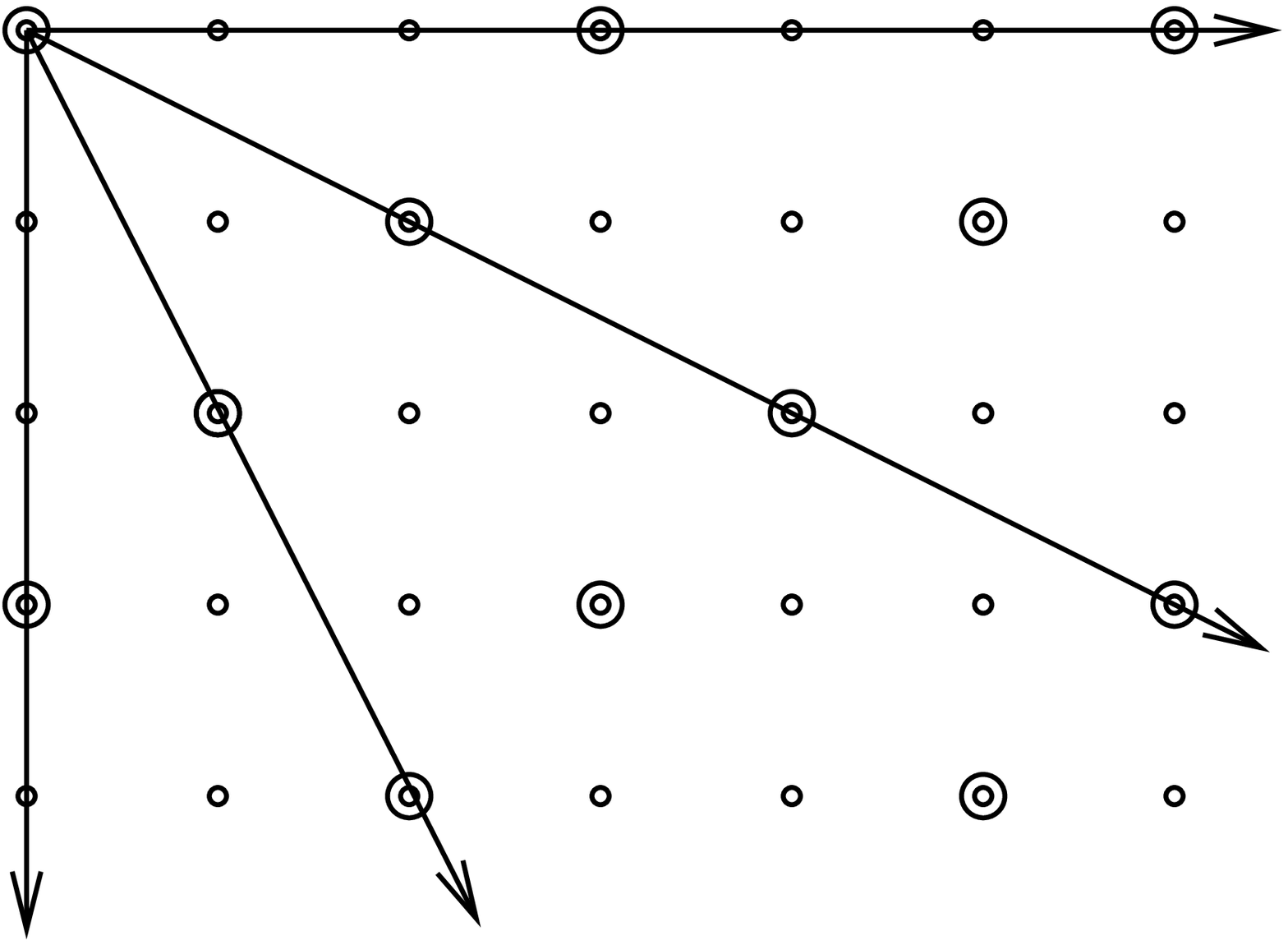}
\end{center}
\caption{\label{Fig1}}
\end{figure}

We claim that locally in $Z$ a toroidal morphism is given by a finite 
toric morphism. A local description of a finite toric morphism can be
seen in Figure \ref{Fig1} where the positive orthant in one lattice 
is subdivided into cones of determinant $1$ in a smaller lattice.
We refer the reader to \cite{Fulton.toric} for the background on
toric geometry.

\begin{remark}\label{subd}
Let $C$ be a positive orthant in a lattice $N$ and let $\hat N$ be 
a finite index sublattice of $N$. Then to each subdivision
$\Sigma$ of $C$ into cones of determinant one in $\hat N$ (see Figure
\ref{Fig1}) one can associate a proper generically finite toric
morphism between the smooth toric variety that corresponds to 
$(\hat N,\Sigma)$ and the smooth toric variety $\CC^{\rk N}$ that 
corresponds to $(N,C)$. 
\end{remark}

As in \cite{Abramovich.Karu}, to every toroidal morphism we associate a map 
$\nu\colon \Sigma_{\hat Z}\to\Sigma_Z$ as follows.
For a cone $C_{\hat I,\hat j}$, pick a generic point $\hat z$ 
on the corresponding connected component $\hat Z_{\hat I;\hat j}$ 
and consider the above monomial map between the neighborhoods of $\hat z$
and $z=\mu (\hat z)$. Point $z$ lies in the stratum $Z_{I;j}$ for some
$I$ and $j$. For every $\hat i\in \hat I$ consider the point in $C_{I;j}$
whose coordinates are the degrees of the local variables that corresponds
to $\hat i$ in the expressions of the variables that correspond to 
$D_i,i\in I$. This defines the image of the point $e_{\hat i;\hat j}$ 
in $C_{I;j}$ and the map is extended to the whole $\Sigma_{\hat Z}$ by 
linearity. This map is compatible with the face inclusions and is thus
well-defined. Indeed, it encodes
the coefficients of $\hat D_{\hat i}$ in the divisors $\mu^*D_i$.

Moreover, we can describe the preimage of any cone $C_{I;j}$ as follows.
Let $z$ be a generic point of the corresponding stratum $Z^\circ_{I;j}$
and let $U$ be a small analytic neighborhood of $z$. Let $\hat U$ be
one of the connected components of the preimage of $U$. We denote by 
$U^\circ$ and $\hat U^\circ$ the intersections of $U$ and $\hat U$
with the complements of $D$ and $\hat D$ respectively. 
Then $\mu$ induces a finite nonramified covering map from $\hat U^\circ$
to $\hat U$. Since the fundamental group of $U^\circ$ is naturally 
isomorphic to $N_{I;j}$, this covering gives a map from $N_{I;j}$ to
a finite group. The kernel of this map is the fundamental group of
$\hat U^\circ$. It is a finite index subgroup of $N=N_{I;j}$
which we denote by $\hat N$. The map from $\hat U$ to $U$ can be passed 
through the singular variety $U_1$ which is the preimage of $U$ under the 
natural map from $\Spec[C_{I;j}^*\cap \hat N]\to \Spec[C_{I;j}^*\cap N]$.
Then the arguments of \cite[Chapter 2,\S 2]{KKMS} 
show that the map from $\hat U$ to
$U_1$ comes from a subdivision of the cone $C_{I;j}$ in $\hat N$ into
cones of volume $1$ as in  Figure \ref{Fig1}.

In general, it is possible that different connected components of the
preimage of $U$ are part of the same connected component of the preimage
of $Z^\circ_{I;j}$. However, we have just shown that the 
preimage $\nu^{-1}C_{I;j}^\circ$ of the interior of $C_{I;j}$
is a union of connected components, each of which corresponds to 
a finite toric morphism. Namely, for each component, there is
a sublattice $\hat N$ of finite index in the lattice $N=N_{I;j}$. 
The relative interior 
of cone $C_{I;j}$ is subdivided into several simplicial cones of volume $1$ 
in the lattice $\hat N$ as in Figure \ref{Fig1}. We will denote this 
subdivision by $\Sigma_{C_{I;j}}$, if it is clear from the 
discussion which connected component of $\nu^{-1}C_{I;j}^\circ$ 
we are referring to. Then every cone of $\Sigma_{C_{I;j}}$ is a cone of 
$\Sigma_{\hat Z}$ with the lattice $\hat N$ as the corresponding lattice.
For the stratum $Z^\circ_{I;j}$ the restriction of $\mu$ to 
the preimage of an analytic 
neighborhood $U$ of $Z^\circ_{I;j}$ is described by 
$\nu$ in the following sense.
Connected components of this preimage are
in one-to-one correspondence with connected components of 
$\nu^{-1}C_{I;j}^\circ$. Let us fix one such connected component which 
we will denote by $\hat U$. It is easily seen to be a neighborhood of
the union of the strata $\hat Z^\circ_{\hat I;\hat j}$ for all
$C_{\hat I;\hat j}\in \Sigma_{C_{I;j}}$. Moreover, it is a locally
trivial fibration with fibers isomorphic to the preimage in 
$\PP_{\hat N,\Sigma_{C_{I;j}}}$ of a disc around
the origin under the map of Remark \ref{subd}. The base is isomorphic
to some smooth variety $W$ which is a nonramified cover of degree
$d_{i;j}$ of $Z^\circ_{I;j}$. The map $\mu$ on $U$ is locally on $Z$
isomorphic to a product of the map from Remark \ref{subd} and an identity
map along. Here $W$ can be taken to be any stratum
on $\hat Z$ that corresponds to a maximum-dimensional cone in 
$\Sigma_{C_{I;j}}$. The numbers $d_{I;j}$ depend on the connected 
component $\hat U$ and they will be important in our description 
of the pushforward.

Our goal is to investigate the restriction of $\mu^*$ and $\mu_*$
to the subrings of $A^*(\hat Z)$ and $A^*(Z)$ which are images of 
$\hat\rho=\rho_{\hat Z}$ and $\rho=\rho_Z$.

\begin{proposition}\label{nustar}
The data of finite toroidal morphism define a map 
$$
\nu^*\colon\CC[[\Sigma_{Z}]]\to\CC[[\Sigma_{\hat Z}]]
$$
simply by pulling back the corresponding functions via $\nu$.
This map is a lifting of $\mu^*$ in the sense that
$$
\hat\rho\circ\nu^*=\mu^*\circ\rho.
$$
\end{proposition}

\begin{proof}
Let's first check this for a linear function on $\Sigma_Z$. Every
such function corresponds to a divisor $\sum_i\alpha_i D_i$.
Then locally this is a statement of toric geometry, namely that
the preimage of a divisor under a map between two toric varieties is 
given by the same piecewise-linear function, which is apparent
from the definition of $\nu$.

It is then enough to check this statement for a function $f$ that 
corresponds to a point $v=\sum_{i\in I}e_{i;j}$ 
in the relative interior of a cone $C_{I;j}\in \Sigma_Z$, since 
both $\nu^*$ and $\mu^*$ are ring homomorphisms. 
We can restrict our attention to a Zariski neighborhood 
$U$ of  $Z_{I;j}$ such that $U\cap Z_I=Z_{I;j}$. Then
all we need is the above statement in toric geometry, restricted to $U$,
together with the fact that $\mu^*$ and $\nu^*$ are ring homomorphisms.
\end{proof}

It is a bit more difficult to describe the pushforward.

\begin{theorem}\label{push}
Let $\nu_*$ be defined as follows. For every $f\in\CC[[\Sigma_{\hat Z}]]$
and every cone $C_{I;j}\in\Sigma_Z$ consider the subdivision $\Sigma_{C_{I;j}}$
of each connected component of $\nu^{-1}C_{I;j}^\circ$. For each 
$C_{\hat I;\hat j}$ of $\Sigma_{C_{I;j}}$ with $\vert \hat I\vert =\vert I\vert $, 
the power series in the variables
$x_{\hat i},\hat i\in \hat I$ that corresponds to the restriction of $f$ 
to $C_{\hat I,\hat j}$ gives a power series in the variables 
$x_i,i\in I$ via the linear change of variables that corresponds to
the inclusion of $C_{\hat I;\hat j}$ into $C_{I;j}$. This power series 
will be denoted by $f_{\hat I;\hat j}$. 
Then we define the element of $\CC((x_{i;j})),i\in I$
\begin{equation}\label{quot}
(\nu_*f)_{I;j}=
\sum_{\Sigma_{C_{I;j}}}d_{I;j}
\sum_{C_{\hat I;\hat j}\in \Sigma_{C_{I;j}},\vert \hat I\vert =\vert I\vert }
f_{\hat I;\hat j}
\,
\frac 
{\prod_{i\in I} x_{i;j}}
{\prod_{\hat i\in \hat I} x_{\hat i;\hat j}}
\end{equation}
where the outer sum is taken over all connected components of 
$\nu^{-1}C_{I;j}^\circ$.
These functions $(\nu_* f)_{I;j}$ are compatible with 
face restrictions and actually lie in $\CC[[x_{i;j},i\in I]]$. 
Thus they define a map $\nu_*\colon\CC[[\Sigma_{\hat Z}]]\to\CC[[\Sigma_Z]]$.
Moreover, $\nu_*$ is a lift of $\mu_*$ in the sense that
$$
\mu_*\circ\hat\rho= \rho\circ\nu_*.
$$
\end{theorem}

\begin{proof}
Let's check that the image of the function $f$ is actually in
$\CC[[C_{I;j}]]$ for every given cone $C_{I;j}$. 
It is enough to consider one cone subdivision $\Sigma_{C_{I;j}}$. The 
function $\nu_* f$ may have simple poles over hyperplanes that correspond to 
cones of $\Sigma_{C_{I;j}}$ of dimension $\vert I\vert -1$. Each such cone has two
adjacent cones of maximum dimension, and it is straightforward to see 
that the terms of \eqref{quot} for two such cones will contribute the 
opposite residues for this hyperplane, because of the compatibility condition
on $f$. 

Let's now show that $\nu_*f$ is well-defined as a map from
$\CC[[\Sigma_{\hat Z}]]$ to $\CC[[\Sigma_Z]]$, that is the definition of
$(\nu_* f)_{I;j}$ is compatible with face inclusions. 
Let $C_{I_1;j_1}$ be a codimension one face of $C_{I_2;j_2}$. 
Hence $I_2=I_1\cup\{i_0\}$.
The cones $C_{\hat I_2;\hat j_2}\in\Sigma_{\hat Z}$ that map to 
the cone $C_{I_2;j_2}$ and have the same dimension may or may not contain 
a face $C_{\hat I_1;\hat j_1}$ that maps to $C_{I_1;j_1}$. 
In the latter case, the contribution of such $C_{\hat I_2,\hat j_2}$
to $(\nu_*f)_{I_2;j_2}$ is going to restrict to zero on the face
$C_{I_1;j_1}$. Indeed, in the equation \eqref{quot}, the restriction 
of the numerator to $C_{I_1;j_1}$ vanishes, because it contains
the linear factor $x_{i_0;j_2}$, whereas the denominator does not vanish. 
In the former case, the factor $x_{i_0;j_2}$ will appear in the numerator
while $x_{\hat i_0;\hat j_2}$ will appear in the denominator where
$\hat i_0$ is the linear function that vanishes on $C_{\hat I_1;\hat j_1}$.
It is easy to see that
$$
\frac{x_{i_0;j_2}}{x_{\hat i_0;\hat j_2}}=
\frac 
{\vert N_{I_2;j_2}:N_{\hat I_2;\hat j_2}\vert }
{\vert N_{I_1;j_1}:N_{\hat I_1;\hat j_1}\vert }.
$$
The other factors in the fraction would restrict
to those for the contribution of $C_{\hat I_1;\hat j_1}$ to
$(\nu_*f)_{I_1;j_1}$. For each cone $C_{\hat I_1;\hat j_1}$ there may be
several different cones $C_{\hat I_2;\hat j_2}$ as above, but we observe
that for each $C_{\hat I_1;\hat j_1}$ 
\begin{equation}\label{deg}
\sum_{C_{\hat I_2;\hat j_2}}{\vert N_{I_2;j_2}:N_{\hat I_2;\hat j_2}\vert }
d_{I_2;j_2}={\vert N_{I_1;j_1}:N_{\hat I_1;\hat j_1}\vert }d_{I_1;j_1}.
\end{equation}
Indeed, both sides depend only on the connected component of 
$\nu^{-1}C_{I_1;j_1}$ rather than the specific cone $C_{\hat I_1;\hat j_1}$.
The right hand side then describes the number of points in the preimage
of a point in the neighborhood of the stratum $Z^\circ_{I_1;j_1}$ 
that lie near a certain connected component 
$Y$ of $\mu^{-1}Z_{I_1;j_1}$. Indeed, for any $z\in Z^\circ_{I_1;j_1}$
there is a neighborhood $U\ni z$ such that the preimage is isomorphic to
a union of $d_{I_1;j_1}$ copies of a toric variety $\PP_{N_{I_1;j_1},
\Sigma_{C_{I_1;j_1}}}$. For each copy the map on the big open set is 
finite unramified of degree 
${\vert N_{I_1;j_1}:N_{\hat I_1;\hat j_1}\vert }$.
 The left hand side describes the sum of numbers of points in 
the preimage of a point in the neighborhood of $Z^\circ_{I_2;j_2}$
sorted by the connected component of its preimage. However, since we are 
summing over the connected components of the preimage of $Z_{I_2;j_2}$ that
are part of the component of the preimage of $Z_{I_1;j_1}$ both sides are 
the same.

Next, we observe that $\nu_*$ is a module homomorphism with respect
to the $\CC[[\Sigma_Z]]$-algebra structure on $\CC[[\Sigma_{\hat Z}]]$
induced by $\nu^*$. Indeed, multiplication by a pullback of a function $g$ on
$\Sigma_Z$ results in multiplication of all $f_{\hat I,\hat j}$ 
in the equation \eqref{quot} by $g$.

Because $\mu_*$ is also a module homomorphism, it is now enough to check 
that $\rho(\nu_*(f))=\mu_*(\hat\rho (f))$ for functions $f$ from some
generating set of $\CC[\Sigma_{\hat Z}]$ as a module over $\CC[\Sigma_Z]$.
We claim that such generating set can be taken to be the set of 
$f$ that correspond to \emph{minimal} lattice points in the interiors of 
cones. Indeed, for every cone $C_{\hat I;\hat j}\in \Sigma_{\hat Z}$, 
consider a cone $C_{I;j}\in\Sigma_Z$ that it maps into the interior of
under $\nu$. It is easy to see that products of the function $f$ 
that corresponds to the minimum interior point of $v=\sum_{\hat i\in I}
e_{\hat i;\hat j}$  by pullbacks of polynomial functions on
$C_{I;j}$ span precisely the space of functions that correspond to 
monomials from the interior of $C_{\hat I;\hat j}$.
So it is now enough to consider the function $f$ that comes from the 
minimum interior point, where we keep the notations as above. We recall
that for such $f$ we have $\hat\rho(f)=\hat Z_{\hat I;\hat j}$.

In the case of $\vert \hat I\vert <\vert I\vert $ the stratum 
$\hat Z_{\hat I,\hat j}$ 
has the image of smaller dimension, so $\mu_*\hat\rho(f)=0$. 
Consider all cones $C_{\hat I_1;\hat j_1}$ in
$\Sigma_{C_{I;j}}$ that contain $C_{\hat I;\hat j}$ and have dimension
$\vert I_1\vert =\vert I\vert $. These cones form the star 
of the neighborhood of 
$C_{\hat I;\hat j}$ in the fan $\Sigma_{C_{I;j}}$. The contributions
of the fractions of the equation \eqref{quot} times the corresponding $f$
will be proportional to
$$
\sum_{C_{\hat I_1;\hat j_1}} \prod_{\hat i\in \hat I_1,i\not\in \hat I}
\frac {1}
{x_{\hat i;\hat j_1}}.
$$
Since all these linear functions $x_{\hat i;\hat j_1}$ vanish on the span
of $C_{\hat I;\hat j}$, one can work on the quotient, where 
Lemma \ref{toricsum} implies that $(\nu_* f)_{I;j}=0$. 

We also need to check
that $(\nu_* f)_{I_1;j_1}=0$ for all cones $C_{I_1;j_1}\in\Sigma_Z$ 
that contain $C_{I;j}$. These calculations are analogous and are left
to the reader.

Let us now consider the case $\vert \hat I\vert =\vert I\vert $. In
this case there will be only one contribution to 
$(\nu_* f)_{I;j}$ and we get
$$
(\nu_* f)_{I;j}=d_{I;j}\prod_{i\in I}x_{i;j},
$$
We also claim that for every other cone $C_{I_1;j_1}$ that contains
$C_{I;j}$ we will have 
$$(\nu_* f)_{I_1;j_1}=d_{I;j}\prod_{i\in I}x_{i;j_1}.
$$
Indeed, for every connected component of $\nu^{-1}(C_{I;j})$ the terms 
of the equation \eqref{quot} will be zero except for the cones 
$C_{\hat I_1;\hat j_1}$ with $\vert \hat I_1\vert =\vert I_1\vert $ that contain
$C_{\hat I,\hat j}$. If we again work in the quotient by the span
of $C_{\hat I,\hat j}$, we see that Lemma \ref{firstorth} implies 
that the total contribution of such $C_{\hat I_1;\hat j_1}$ is
$$
d_{I_1;j_1}
\frac
{\vert N_{I_1;j_1}:N_{\hat I_1;\hat j_1}\vert }
{\vert N_{I;j}:N_{\hat I;\hat j}\vert }
\prod_{i\in I}x_{i;j_1}.
$$
Then an analog of equation \eqref{deg} finishes the calculation of 
$(\nu_* f)_{I_1;j_1}$.

Then we conclude that $\nu_*f$ corresponds to $d_{I;j}$ times the 
minimum lattice point of $C_{I;j}$. On the other hand,
the corresponding cycle $\hat Z_{\hat I,\hat j}$ maps onto 
$Z_{I;j}$ and the morphism is is generically finite of degree
$d_{I;j}$. Therefore, we get 
$\mu_*\hat Z_{\hat I,\hat j}=d_{I;j}Z_{I;j}$, which finishes the proof.
\end{proof}

\section{Main theorem}\label{section.main}
Our goal is to prove that any $G$-equivariant morphism $\mu:\hat Z\to Z$ 
of smooth varieties which is birational to a quotient by $G$ 
has the property that the pushforward of the orbifold elliptic 
class in $A^*(\hat Z)$ is the elliptic class in $A^*(Z)$. The strategy 
of the proof is to reduce the situation to a toroidal morphism.

Let $\mu:\hat Z\to Z$ be such  $G$-equivariant 
toroidal morphism. Let $h$ be any linear function
on the fan $\Sigma_Z$. Function $h$ corresponds to the divisor
$$
D_h=\sum_{i\in I_Z} \alpha_iD_i.
$$
We will also consider the divisor $\hat D_h$ on $\hat Z$ 
such that 
\begin{equation}\label{discr}
\mu^*(K_Z+D_h)=K_{\hat Z}+\hat D_h
\end{equation}
and 
$$
\hat D_h=\sum_{\hat i\in I_{\hat Z}} \hat \alpha_{\hat i} \hat D_{\hat i}.
$$
The set $I_{\hat Z}$ splits into two sets $I_{\hat Z}^{exc}$ and 
$I_{\hat Z}^{ram}$ according to whether or not $D_i$ is contracted by
$\mu$ to a smaller-dimensional variety. We will abuse notations and call
the divisors $\hat D_{\hat i}$ for $i\in I_{\hat Z}^{ram}$ ramification 
divisors, and assign to each of them the ramification index $r_{\hat i}$ 
(which may be equal to $1$). We observe that if 
$\mu(\hat D_{\hat i})=D_i$, then
$$
\hat \alpha_{\hat i}+1 = r_{\hat i}(\alpha_i +1).
$$
Since $\mu$ is birationally equivalent to a quotient morphism and 
is locally given by the map between two toric varieties, it is easy
to see that locally the group $G$ acts as a subgroup of the torus.
The isotropy group of every point of the stratum that corresponds to 
the cone $C_{\hat I;\hat j}$ is equal to the index of 
$\sum_{\hat i\in \hat I}\ZZ e_{\hat i;\hat j}$ in the restriction 
of $N_{C_{I;j}}$ to the $\QQ$-span of $e_{\hat i;\hat j}$.
We will denote this group by $G_{\hat I;\hat j}$. 

Consider the following function $E$ of $z,\tau$ with values in $A^*(\hat Z)$.
$$
E(z,\tau)=
\frac 1{\vert G\vert } \sum_{g,h\in G;gh=hg} \sum_{Z^{g,h}=Z_{\hat I;\hat j}} 
Z_{\hat I;\hat j}
~
\prod_{\hat i\in I_{\hat Z}} 
\frac
{\hat D_{\hat i}\theta(\frac{\hat D_{\hat i}}{2\pi\ii}-z)\theta'(0)}
{2\pi\ii\theta(\frac{\hat D_{\hat i}}{2\pi\ii})\theta(-z)}
$$
$$
\prod_{\hat i\in \hat I}
\frac
{\theta(\frac{\hat D_{\hat i}}{2\pi\ii}+g_{\hat i}-h_{\hat i}\tau-
(\hat\alpha_{\hat i}+1)z)
\theta(-z)
\theta(\frac{\hat D_{\hat i}}{2\pi\ii})
}
{\hat D_{\hat i}
\theta(\frac{\hat D_{\hat i}}{2\pi\ii}+g_{\hat i}-h_{\hat i}\tau)
\theta(-(\hat \alpha_{\hat i}+1)z)
\theta(\frac{\hat D_{\hat i}}{2\pi\ii}-z)}
\ee^{2\pi\ii (\hat \alpha_{\hat i}+1)h_{\hat i} z}
$$
$$
\prod_{\hat i\not\in\hat I}
\frac
{\theta(\frac{\hat D_{\hat i}}{2\pi\ii}-(\hat\alpha_{\hat i}+1)z)
\theta(-z)}
{\theta(\frac{\hat D_{\hat i}}{2\pi\ii}-z)\theta(-(\hat\alpha_{\hat i}+1)z)}.
$$
Here $g_{\hat i}$ and $h_{\hat i}$ are the rational numbers in the range
$[0,1)$ which describe the characters of the action of $g$ and $h$ on
the divisor $\hat D_{\hat i}$ at each point of $\hat Z_{\hat I;\hat j}$.

The following theorem describes the pushforward of $E$ to $Z$.
\begin{theorem}\label{bigmess}
$$\mu_*E(z,\tau) =
\prod_{i\in I_Z}
\frac
{D_i
\theta'(0)
\theta(\frac{D_{i}}{2\pi\ii}-(\alpha_{i}+1)z)
}
{2\pi\ii\,
\theta(\frac {D_i}{2\pi\ii})
\theta(-(\alpha_{i}+1)z)}.
$$
\end{theorem}

\begin{proof}
We will use Theorem \ref{push} to reduce the statement to 
a combinatorial result. First, we observe that $E(z,\tau)$ 
can be obtained as $\hat\rho(F(z,\tau))$ where $F$ is defined as follows.

Consider the cone $C_{\hat I;\hat j}$ that is a part of the 
subdivision $\Sigma_{C_{I;j}}$. Denote by $G_{\hat I;\hat j}$ 
the quotient of the intersection of the lattice $\oplus_{i\in I}\ZZ e_{i;j}$ 
with the rational span of $\hat e_{\hat i;\hat j},\hat i\in\hat I$ 
by  $\oplus_{\hat i\in\hat I}\ZZ \hat e_{\hat i;\hat j}$.
The value of $F(z,\tau)$ on this cone is 
$$
F_{\hat I;\hat j}(z,\tau)=
\frac 1{\vert G\vert } \sum_{g,h\in G_{\hat I;\hat j}}
\prod_{\hat i\in \hat I} 
\frac
{
\hat x_{\hat i;\hat j}
\theta'(0)
\theta(\frac{\hat x_{\hat i;\hat j}}{2\pi\ii}+g_{\hat i}-h_{\hat i}\tau-
(\hat\alpha_{\hat i}+1)z)
}
{
2\pi\ii\,
\theta(\frac{\hat x_{\hat i;\hat j}}{2\pi\ii}+g_{\hat i}-h_{\hat i}\tau)
\theta(-(\hat \alpha_{\hat i}+1)z)
}
\ee^{2\pi\ii (\hat \alpha_{\hat i}+1)h_{\hat i} z}.
$$
Indeed, for every $g,h\in G_{\hat I;\hat j}$ there is a unique connected
component of the fixed point of $g$ and $h$ that contains $Z_{\hat I;\hat j}$
Then we use the definition and multiplicative properties of $\rho$, together
with the fact that $\hat D_{\hat i}$ corresponds to the function on 
$\Sigma_{\hat Z}$ that equals $x_{\hat i,\hat j}$ for every $C_{\hat I;\hat j}$
such that $\hat I\ni \hat i$ and equals zero otherwise. 

Now we need to calculate $\nu_*F(z,\tau)$. Let's calculate the 
component of $\nu_*F(z,\tau)$ on a cone $C_{I;j}$. 
By the definition of $\nu_*$ we get
$$
\nu_*F(z,\tau)_{I;j}=
({\prod_{i\in I} x_{i;j}})
\sum_{\Sigma_{C_{I;j}}}
\frac {d_{I;j}}{\vert G\vert }
$$
$$
\sum_{C_{\hat I;\hat j}\in \Sigma_{C_{I;j}},\vert \hat I\vert =\vert I\vert }
\sum_{g,h\in G_{\hat I;\hat j}}
\prod_{\hat i\in \hat I} 
\frac
{
\theta'(0)
\theta(\frac{\hat x_{\hat i;\hat j}}{2\pi\ii}+g_{\hat i}-h_{\hat i}\tau-
(\hat\alpha_{\hat i}+1)z)
}
{
2\pi\ii\,
\theta(\frac{\hat x_{\hat i;\hat j}}{2\pi\ii}+g_{\hat i}-h_{\hat i}\tau)
\theta(-(\hat \alpha_{\hat i}+1)z)
}
\ee^{2\pi\ii (\hat \alpha_{\hat i}+1)h_{\hat i} z}.
$$
We now apply Lemma \ref{mainthetalemma}. Indeed, it is easy to check 
that equation \eqref{discr} implies that the values $(\alpha_{\hat i}+1)z$ 
are values of a linear function on $C_{I;j}$.
As a result, we get
$$
\nu_*F(z,\tau)_{I;j}=
({\prod_{i\in I} x_{i;j}})
\sum_{\Sigma_{C_{I;j}}}
\frac {d_{I;j}}{\vert G\vert }
\vert N_{I;j}:N_{\hat I;\hat j}\vert 
\prod_{i\in I} 
\frac
{
\theta'(0)
\theta(\frac{x_{i;j}}{2\pi\ii}-(\alpha_i+1)z)
}
{
2\pi\ii\,
\theta(\frac{x_{i;j}}{2\pi\ii})
\theta(-(\alpha_i+1)z)
}
$$
$$
=\prod_{i\in I}\frac
{x_{i;j}
\theta'(0)
\theta(\frac{x_{i;j}}{2\pi\ii}-(\alpha_i+1)z)
}
{
2\pi\ii\,
\theta(\frac{x_{i;j}}{2\pi\ii})
\theta(-(\alpha_i+1)z)
}.
$$
Here we use $\sum_{\Sigma_{C_{I;j}}}d_{I;j}\vert N_{I;j}:N_{\hat I,\hat j}\vert =\vert G\vert $,
which follows from the count of the number of preimage points of a point close
to the stratum $Z_{I;j}$.

We now use Theorem \ref{push} to get
$$
\mu_*E(z,\tau)=\mu_*\hat\rho F(z,\tau)=\rho\nu_*F(z,\tau)
=
\prod_{i\in I_Z}
\frac
{D_i
\theta'(0)
\theta(\frac{D_{i}}{2\pi\ii}-(\alpha_{i}+1)z)
}
{2\pi\ii\,
\theta(\frac {D_i}{2\pi\ii})
\theta(-(\alpha_{i}+1)z)}.
$$
Indeed, the calculation of $\rho\nu_*F(z,\tau)$ is accomplished by 
the multiplicativity of $\rho$ and the fact that the power series in $D_i$
has constant term $1$.
\end{proof}

We will also need the following lemma that connects the Chern classes
of $T\hat Z$ and $\mu^*TZ$. 
\begin{lemma}\label{log}
$$
c(T\hat Z)=\prod_{\hat i\in I_{\hat Z}}(1+D_{\hat i})
\prod_{i\in I}(1+\mu^*D_i)^{-1}\mu^*c(TZ).
$$
\end{lemma}

\begin{proof}
First of all, it is easy to see that the pullback of the bundle of the 
logarithmic differentials on $Z$ is the bundle of logarithmic
differentials on $\hat Z$. Then it is straightforward to calculate the 
ratio of Chern classes for the bundles of logarithmic differentials 
and usual differentials for a variety with the normal crossing divisor. 
The details are left to the reader.
\end{proof}

We are now ready to formulate and prove our main theorem.
\begin{theorem}\label{main}
Let $(X;D_X)$ be a Kawamata log-terminal pair which is invariant under 
an effective action of $G$ on $X$. Let 
$
\psi\colon X\to X/G 
$
be the quotient morphism. Let  $(X/G;D_{X/G})$ be the quotient pair,
see Definition \ref{quotpair}.
Then 
$$
\psi_* {\cal Ell}_{orb}(X,D_X,G;z,\tau)={\cal Ell}(X/G,D_{X/G};z,\tau).
$$
\end{theorem}

\begin{proof}
The following lemma allows us to reduce the problem to the situation
of a $G$-equivariant toroidal morphism.

\begin{lemma}\label{reduce}
There exists a commutative diagram
$$
\begin{array}{rccc}
\mu\colon&\hat Z&\to&Z\\
 &\downarrow& &\downarrow\\
\psi\colon &X&\to& X/G
\end{array}
$$
where the vertical arrows are resolutions of singularities and 
$\mu$ is a $G$-equivariant toroidal morphism.
\end{lemma}

\begin{proof}
We define $Z$ as a desingularization of $(X/G,D_{X/G})$. Consider 
the normalization of $Z$ in the function field of $X$ and the
corresponding normalization morphism. By Abhyankar's
lemma it is a (typically singular) toroidal embedding with the
toroidal morphism to $Z$. Then toroidal desingularization finishes
the job. See \cite{Abramovich.Wang} for details.
\end{proof}

\emph{Proof of Theorem \ref{main} continues.}
By Lemma \ref{reduce}, Definition \ref{genorb} and composition
properties of pushforwards, it is 
sufficient to prove the pushforward 
result for a $G$-equivariant toroidal morphism
$\mu\colon \hat Z\to Z$ which is birational to $\psi$. 
By Lemma \ref{log} and Definition \ref{maindef}
of orbifold elliptic class ${\cal Ell}
(\hat Z,D_{\hat Z},G;z,\tau)$, we see that 
$$
{\cal Ell}(\hat Z,D_{\hat Z},G;z,\tau)=E(z,\tau)\mu^*
\Big(
\prod_k
\frac{z_k\theta(\frac {z_k}{2\pi\ii}-z)}{\theta(\frac{z_k}{2\pi\ii})}
\prod_{i\in I_Z}
\frac
{2\pi\ii\theta(\frac{D_i}{2\pi\ii})\theta(-z)}
{D_i\theta'(0)\theta(\frac{D_i}{2\pi\ii}-z)}
\Big).
$$
Then Theorem \ref{bigmess} and definition of elliptic class
${\cal Ell}(Z,D;z,\tau)$ finishes the proof.
\end{proof}

\begin{remark}
Theorem \ref{main} gives an affirmative answer to the conjecture of \cite
{singellgenus}. We call it McKay correspondence for elliptic genera,
analogously to the homological McKay correspondence for stringy $E$-functions.
\end{remark}

\section{DMVV formula for pairs}\label{section.DMVV}
One of the motivations of the definition of orbifold elliptic genus in
\cite{singellgenus} was the formula for the generating functions of elliptic
genera of symmetric products. 
\begin{equation}
\sum_{n\geq 0}p^n {Ell}_{orb}(X^n/S_n;z,\tau)=
\prod_{i=1}^\infty \prod_{l,m}
\frac 1
{(1-p^iy^lq^m)^{c(mi,l)}}.
\end{equation}
Here $X$ is a K\"{a}hler manifold, 
$S_n$ is the symmetric group acting on the
$n$-fold product and $c(m,l)$  are the coefficients of 
the elliptic genus $\sum_{m,l}c(m,l)y^lq^m$ of $X$. 

This formula was originally derived in \cite{DMVV}
by means of some string-theoretic arguments. In particular, the orbifold
elliptic genus of a quotient of a variety $X^n$ by the symmetric group
$S_n$ was defined as trace of a certain operator 
over the Hilbert space of the
conformal field theory quotient of ${\cal C}^n$ where $\cal C$ is 
the superconformal field theory conjecturally associated to $X$. 
In \cite{singellgenus} DMVV formula was shown for the mathematically defined
orbifold elliptic genus. Our goal now is to extend this result to
singular varieties and more generally to arbitrary Kawamata log-terminal pairs.

\begin{theorem}\label{DMVVpairs}
Let $(X,D)$ be a Kawamata log-terminal pair. For every $n\geq 0$ consider
the quotient of $(X,D)^n$ by the symmetric group $S_n$, which 
we will denote by $(X^n/S_n,D^{(n)}/S_n)$. Here we denote by 
$D^{(n)}$ the sum of pullbacks of $D$ under $n$ canonical projections to $X$. 
Then we have 
$$
\sum_{n\geq 0}p^n {Ell}(X^n/S_n,D^{(n)}/S_n;z,\tau)=
\prod_{i=1}^\infty \prod_{l,m}
\frac 1 {(1-p^iy^lq^m)^{c(mi,l)}},
$$
where the elliptic genus of $(X,D)$ is
$$
\sum_{m\geq 0}\sum_l c(m,l)y^lq^m
$$
and $y=\ee^{2\pi\ii z}$, $q=\ee^{2\pi\ii\tau}$.
\end{theorem}

\begin{remark}
In the case of smooth $X$ with $D=0$, the Fourier coefficient at of
$Ell(X,z,\tau)$ at $q^m$ is a polynomial in $y^{\pm \frac12}$.
In general other rational $l$ are possible, but more importantly, 
the coefficient by $q^m$ is no longer a polynomial in 
$y^{\pm \frac1d}$, rather it is a rational function. However, we will
always assume that this function is Laurent expanded around $y=0$,
so we will be working in the field of formal power series in 
$y^{\pm \frac 1d}$, where $d$ is divisible by $2$ and all denominators
of the discrepancy coefficients for some resolution of $(X,D)$.
The issue of non-polynomiality was first raised in
\cite{Batyrev.1} at the $q^0$ level.
\end{remark}

\begin{proof}
First of all, observe that the quotient of the tensor power of a 
Kawamata log-terminal pair is again a Kawamata log-terminal pair.
Moreover, by Theorem \ref{main}, we can calculate the elliptic genus
of $(X^n/S_n,D^{(n)}/S_n)$ as an orbifold elliptic genus of 
$(X^n,D^{(n)},S_n)$. 
A resolution of singularities $(\hat X,\hat D)$ of the pair $(X,D)$
induces a birational morphism $\hat X^n\to X^n$ so we may assume 
that $(X,D)$ is nonsingular, i.e. $X$ is smooth and $D$ is a normal crossing
divisor $\sum_i \alpha_iD_i$ with $\alpha_i> -1$.
While the divisor $D^{(n)}$ on $X^n$ has simple normal crossings, it
is not $S_n$-normal. Indeed the pullbacks of the same
component $D_i$ via different projections are group translates of each other
and certainly intersect and are nontrivially permuted by the isotropy group of 
any such intersection point that lies on the main diagonal $X\subseteq X^n$.
To rectify this situation we need to consider an appropriate blowup of $X^n$. 
By Remark \ref{eachgh}, each pair of commuting elements
$(g,h)$ can be handled separately. 

Let's describe the pairs of commuting elements $g,h\in S_n$
and the connected components of their fixed point set. If the cycle
decomposition of $h$ has $a_j$ cycles of degree $j$,
then the fixed point set of $h$ on $X^n$ is the product of 
$\sum_j a_j$ copies of $X$, embedded into $X^n$ by the product of 
diagonal embeddings of $X$ into $X^j$ for each cycle of length $j$.
Elements $g$ of $S_n$ that commute with $h$ form a semidirect product 
of the group $C_h=\prod_j (\ZZ/j\ZZ)^{a_j}$ which consists of the products 
of powers of cycle components of $h$ and the group $B_h=\prod_jS_{a_j}$ which
consists of the group that permutes cycles of the same length without 
disturbing the order in the cycle. A fixed point set of each such pair
$(g,h)$  consists of points on $X^{\sum_j a_j}$ that are preserved by 
the image of $g$ in $B_h$. 
It is easy to see that the contribution of each such $(X^n)^{g,h}$ is 
the product of the contributions of each factor. As a result, it is 
enough to consider the contribution of the diagonal embedding of $X$
into $X^{ij}=(X^j)^i$ where $h$ acts by permuting the copies of 
$X$ inside each $X^j$ and $g$ acts by a product of a cyclic permutation
of $i$ copies of $X^j$ and some cyclic permutations within each $X^j$,
that does not change the cyclic orders of the components of $X^j$.
Then $g^i=h^s$ for some $0\leq s\leq j-1$, and $s$ determines the action
uniquely. Namely, if  $x_{k,l},k\in\ZZ/i\ZZ,l\in\ZZ/j\ZZ$  
denote the components of $X^{ij}$, then we may assume that 
$h$ acts by $x_{k,l}\to x_{k,l+1}$ and $g$ acts by 
$x_{k,l}\to x_{k+1,l}$ for $k=0,\ldots, i-2$ and $x_{i-1,l}\to x_{0,l+s}$. 
We will denote by $G$ the group generated by $g$ and $h$. It is an abelian
group of order $ij$ given by the generators $g,h$ and relations
$gh=hg, g^i=h^s, h^j=1$. We denote the corresponding product
of $ij$ copies of $X$ by $X^G$, which indicates the action of $G$ on it.

We now need to make $(X^{G},D^{(G)})$ into a $G$-normal pair.
Let $D_c,1\leq c\leq k$ be the irreducible components of $D$ on $X$. 
We will denote by $D_{r,c},r\in G$ the pullback of $D_c$ under the 
$r$-th projection map $X^{G}\to X$. We will perform the following sequence 
of blowups to $X^G$. First, we blow up $\cap_{r\in G}D_{r,1}$, then we 
blow up the proper preimage of $\cap_{r\in G}D_{r,2}$, and so on. 
We can describe this blowup in terms of the subdivision of the conical 
complex that corresponds to the simple normal crossing 
divisor $D^{(G)}$ on $X^G$.
For the sake of simplicity we assume that the intersection of every number 
of components $D_c$ on $X$ is connected. The general case is completely
analogous, it can also be reduced to the connected case by further
blowups of $X$. Every cone $C$ of the conical complex $\Sigma_{X^G}$ is 
generated by elements $e_{r,c}$ for some subset of $I_C\subseteq 
G\times \{1,\ldots,k\}$. We denote by $J_C$ the subset of $\{1,\ldots,k\}$ 
that consists of all $c$ for which $I_C\supseteq G\times \{c\}$. The 
subdivision of $C$ that corresponds to this sequence of blowups 
is then the product of $\ZZ_{\geq 0}e_{r,c}$ for $(r,c)\in I_C,c\not\in
J_C$ and the product over all $c\in J_C$ of the subdivisions of
$\sum_{r\in G}\ZZ_{\geq 0}e_{r,c}$ where the extra vertex 
$\sum_{r\in G}e_{r,C}$ is added and the cone is subdivided accordingly.
It is clear that this is a well-defined subdivision of $\Sigma_{X^G}$ 
and we denote the corresponding variety by $\hat X^G$ and the corresponding
divisor by $\hat D^{(G)}$. We observe that there are $k$ exceptional components
of $\hat D^{(G)}$, which we will call $E_c$,
and the rest are the proper preimages of the components of $D^{(G)}$.

We need to describe connected components of the fixed point set of 
$G$ on $\hat X^G$. Every such fixed point maps to the diagonal 
$X\subseteq X^G$, and should lie on the stratum of the stratification by the 
intersections of components of $\hat D^{(G)}$ that is stable under
the group action. Since the construction is local in $X$, we need to see
what happens when $X$ is a $\CC^n$
with $D$ given as a union of some coordinate hyperplanes 
$z_1=0,z_2=0,\ldots, z_l=0$. The extra coordinates $z_{l+1},\ldots, z_n$ 
will have an effect of tensoring the construction
by an affine space, so it is enough to look at $l=n$ case. 
Then we need to investigate the fixed point sets of the toric variety that 
corresponds to a certain blowup of the positive orthant in $\ZZ^{ijk}$ where 
the generators are denoted by $e_{r,c},r\in G,1\leq c\leq l$.
The group $G$ acts by multiplication on the first component of the index of
the coordinate. The rays of the fan of the blowup that are fixed under $G$ 
correspond to $e_{*,c}=\sum_{r\in G}e_{r,c}$. Moreover, it is easy to see
that the only strata that are preserved by $G$ are the intersections of 
the corresponding divisors.
In other words, we need to consider the faces of the $l$-dimensional cone $C$ 
which is a part of the subdivision of the positive orthant and 
is the span of all $e_{*,c}$. This cone corresponds to the affine set
which is isomorphic to
\begin{equation}\label{loc}
\CC^{l}\times (\CC^*)^{ijl-l}.
\end{equation}
The coordinates on $(\CC^*)^{ijl-l}$ are given by 
$
x_{r,c}x_{r_1,c}^{-1}
$
and the coordinates on $\CC^l$ are given by $x_{0,c}$.
Let $P$ be a fixed point of $G$. For each $c$,
$x_{r,c}/x_{r_1,c}=\exp(2\pi\ii\lambda(r-r_1))$ for some character 
$\lambda\colon G\to\QQ/\ZZ$.
If $\lambda$ is nontrivial then $x_{0,c}$ is zero, and otherwise
arbitrary values of $x_{0,c}$ are allowed. Moreover, for each component 
of the fixed point set the map to $\CC^n\subseteq (\CC^n)^G$ is an embedding. 
Indeed, it is clear for each factor $\CC^G$ that corresponds to the $D_c$.
Basically, for each factor, the blowup locus intersects the main diagonal
of $\CC^G$ in codimension one, namely at the origin.

Returning to the global situation, the above description tells us that
connected components of the fixed point set $Y$ of $X^{(G)}$ correspond to 
the collections of characters $\lambda_c\colon G\to \QQ/\ZZ$.
The fixed point set for each such character is isomorphic to
$D_I=\cap_{i\in I}D_i$ where $I$ is the set of $c$-s 
for which $\lambda_c$  is nontrivial. Indeed, this follows from the 
fact that locally the map from the component of the fixed point set to 
$X^G$ is an embedding. We observe that for some combinations
of characters we may have $D_I=\emptyset$. 

We now need to calculate the tangent bundle to such a component
which we will denote by $Y_{\lambda_1,\ldots,\lambda_k}$.
Notice that divisors $\hat D_{r,c}$ do not intersect with 
$Y$. Indeed, every $G$-invariant point of $\hat D_{r,c}$ would belong
to $\hat D_{r_1,C}$ for all $r_1$, but the intersection of all these
divisors is empty since $\pi$ factors through the blowup of the intersection
of $D_{r,C},r\in G$. As far as intersection with $E_c$ is concerned, 
$Y_{\lambda_1,\ldots,\lambda_k}$ is contained in $E_c$ for $\lambda_c\neq 0$
and intersects transversally other $E_c$. For $\lambda_c=0$ the intersection
of $E_c$ and $Y_{\lambda_1,\ldots,\lambda_k}$ can be identified with 
the intersection by $D_c$ under the isomorphism 
$Y_{\lambda_1,\ldots,\lambda_k}\cong D_I$. The character of $G$ that 
corresponds to $E_c\supseteq Y_{\lambda_1,\ldots,\lambda_k}$
is equal to $\lambda_c$. 

The Chern classes of the tangent bundles of $\hat X^G$ and $X^G$ 
are related by Lemma \ref{log}, namely 
$$
c(T\hat X^G)=\pi^*c(TX^G)\prod_{c=1}^k(1+E_c)\prod_{r\in G,1\leq c\leq k}
\frac {(1+\hat D_{r,c})}{(1+\pi^*D_{r,c})}
$$
where $\hat D_{r,c}$ is the proper preimage of $D_{r,c}$. Notice
that as classes in $A^*(\hat X^G)$, $\hat D_{r,c}=\pi^* D_{r,c}-E_c$.
Moreover, we can write $c(TX^G)$ as $\oplus_{r\in G} TX_r$ where 
$TX_r$ is the pullback of the tangent bundle of $X$ under the $r$-th
projection. Since $\hat D_{r,c}$ are disjoint from 
$Y_{\lambda_1\ldots,\lambda_k}$, we get
$$
c(i^*T\hat X^G)=i^*\pi^*c(TX^G)
\prod_{c=1}^k(1+i^*E_c)^{1-\vert G\vert }
$$
where $i\colon Y_{\lambda_1,\ldots,\lambda_k}\to \hat X^{G}$ is the 
embedding. Notice that $\pi$ restricts to an embedding on 
$Y_{\lambda_1,\ldots,\lambda_k}$ with image $D_I\subseteq X\subseteq X^G$
where $I$ is the set of all $c$ that for which $\lambda_c$ is nontrivial.
The following lemma describes $i^*T\hat X^G$ in more detail.

\begin{lemma}\label{tantoy}
Let $\lambda$ be a character of $G$. Then the $\lambda$-component 
$V_\lambda$ of 
the restriction of $T\hat X^G$ to $Y_{\lambda_1,\ldots,\lambda_k}$,
identified with $D_I\neq \emptyset$ can be described as follows.
If $\lambda=0$, then $V_\lambda=TD_I$. If $\lambda\neq 0$, then
there is an exact sequence
$$0\to j^*TX_{\log}\to V_\lambda\to
\bigoplus_{c,\lambda_c=\lambda}{\cal O}(D_c)\to 0$$
where $j$ is the embedding $D_I\to X$ and $TX_{\log}$ is the dual to 
the bundle of log-differentials for $(X,D)$.
\end{lemma}

\begin{proof}
We observe that $Y_{\lambda_1,\ldots,\lambda_k}$ 
is contained in the intersection of
$E_c$ for $\lambda_c\neq 0$, which induces a $G$-equivariant surjection from
the restriction of $T\hat X^G$ to the restriction of
$\oplus_{\lambda_c\neq 0} {\cal O}(E_c)$ with the kernel being the 
restriction of the tangent space to $\cap_{\lambda_c\neq 0}E_c$ to
$Y_{\lambda_1,\ldots,\lambda_k}$. It is 
easy to see that under the identification of
$Y_{\lambda_1,\ldots,\lambda_k}$ with $D_I$ the restriction of ${\cal O}E_c$ 
is isomorphic to ${\cal O}(D_c)$ and has character $\lambda_c$.

So we now need to investigate the restriction of the tangent space of
$\cap_{\lambda_c\neq 0}E_c$ and its eigenbundles.
The $\lambda=0$ case is clear, so it is enough to consider the normal
bundle to $Y_{\lambda_1,\ldots,\lambda_k}$ in $\cap_{\lambda_c\neq 0}E_c$.
Locally in the notations of \eqref{loc}, this bundle is isomorphic
to the restriction of
the tangent bundle of $(\CC^*)^{ijk-k}$. The cotangent bundle of 
$(\CC^*)^{ijk-k}$
is generated by 
$\frac {dx_{r,c}}{x_{r,c}}-\frac {dx_{r_1,c}}{x_{r_1,c}}$, so its $\lambda$-eigenbundle
is isomorphic to a bundle generated by $\frac {dx_{c}}{x_c}$, which is precisely
the bundle of logarithmic differential forms. Even though \eqref{loc} refers to the neighborhood
of a point of the intersection of $\dim X$ divisors $D_c$, it is clear that the general case
is obtained by a Cartesian product with a disc and the identification is still valid.
It remains to notice that this identification behaves well under coordinate changes. 
\end{proof}

\emph{Proof of Theorem \ref{DMVVpairs} continues.}
In view of Lemma \ref{tantoy}, the contribution of $(g,h)$ to the 
orbifold elliptic genus of $(X^G,D^{(G)})$ is 
$$
\sum_{\{\lambda_1,\ldots,\lambda_k\},\cap_{\lambda_c\neq 0}D_c\neq\emptyset}
\int_X 
\Big(\prod_{c,\lambda_c\neq 0} D_c\Big)
\prod_{l}\frac {x_l\theta(\frac {x_l}{2\pi\ii}-z)}
{\theta(\frac {x_l}{2\pi\ii})}
\prod_{c,\lambda_c\neq 0}
\frac
{\theta(\frac {D_c}{2\pi\ii})}
{D_c\theta(\frac{D_c}{2\pi\ii}-z)}
$$
$$
\times
\prod_{\lambda\neq 0}
\Big(\prod_{l}
\frac {\theta(\frac {x_l}{2\pi\ii}+\lambda(g)-\lambda(h)\tau-z)}
{\theta(\frac {x_l}{2\pi\ii}+\lambda(g)-\lambda(h)\tau)}
\ee^{2\pi\ii\lambda(h)z}
$$
$$\times
\prod_{c=1}^k
\frac
{\theta(\frac{D_c}{2\pi\ii}+\lambda(g)-\lambda(h)\tau)
\theta(\lambda(g)-\lambda(h)\tau-z)}
{\theta(\frac{D_c}{2\pi\ii}+\lambda(g)-\lambda(h)\tau-z)
\theta(\lambda(g)-\lambda(h)\tau)}
\Big)
$$
$$
\times
\prod_{c,\lambda_c\neq 0}
\frac {\theta(\frac {D_c}{2\pi\ii}+\lambda_c(g)-\lambda_c(h)\tau-z)}
{\theta(\frac {D_c}{2\pi\ii}+\lambda_c(g)-\lambda_c(h)\tau)}
\ee^{2\pi\ii\lambda_c(h)z}
$$
$$\times
\prod_{c,\lambda_c\neq 0}
\frac 
{\theta(\frac {D_c}{2\pi\ii}+\lambda_c(g)-\lambda_c(h)\tau-\vert G\vert (\alpha_c+1)z)
\theta(-z)}
{\theta(\frac {D_c}{2\pi\ii}+\lambda_c(g)-\lambda_c(h)\tau-z)
\theta(-\vert G\vert (\alpha_c+1)z)}
\ee^{2\pi\ii(\vert G\vert \alpha_c+\vert G\vert -1)\lambda_c(h)z}
$$
$$\times
\prod_{c,\lambda_c=0}
\frac 
{\theta(\frac {D_c}{2\pi\ii}-\vert G\vert (\alpha_c+1)z)
\theta(-z)}
{\theta(\frac {D_c}{2\pi\ii}-z)
\theta(-\vert G\vert (\alpha_c+1)z)}
$$
where we have used the fact that the coefficients by $E_c$ in the 
log-pair on $\hat X^G$ are $(\vert G\vert \alpha_c-\vert G\vert -1)$ and other divisors
do not intersect the fixed point set and are thus irrelevant.
After observing that the formula gives $0$ for the case 
$D_I=\emptyset$, the above can be rewritten as
$$
F_{i,j,s}=\sum_{\{\lambda_1,\ldots,\lambda_k\}}
\int_X \prod_l \Big(x_l\prod_{\lambda}
\frac {\theta(\frac {x_l}{2\pi\ii}+\lambda(g)-\lambda(h)\tau-z)}
{\theta(\frac {x_l}{2\pi\ii}+\lambda(g)-\lambda(h)\tau)}
\ee^{2\pi\ii\lambda(h)z}
\Big)
$$
$$
\times
\prod_{\lambda\neq 0}
\prod_{c=1}^k
\frac
{\theta(\frac{D_c}{2\pi\ii}+\lambda(g)-\lambda(h)\tau)
\theta(\lambda(g)-\lambda(h)\tau-z)}
{\theta(\frac{D_c}{2\pi\ii}+\lambda(g)-\lambda(h)\tau-z)
\theta(\lambda(g)-\lambda(h)\tau)}
$$
$$
\times
\prod_{c}
\frac 
{\theta(\frac {D_c}{2\pi\ii})
\theta(\frac {D_c}{2\pi\ii}+\lambda_c(g)-\lambda_c(h)\tau-\vert G\vert (\alpha_c+1)z)
\theta(-z)}
{\theta(\frac{D_c}{2\pi\ii}-z)
\theta(\frac {D_c}{2\pi\ii}+\lambda_c(g)-\lambda_c(h)\tau)
\theta(-\vert G\vert (\alpha_c+1)z)}
\ee^{2\pi\ii\vert G\vert (\alpha_c+1)\lambda_c(h)z}.
$$

We will use the following lemmas that take into account the specific 
form of $G$.
\begin{lemma}\label{product}
$$\prod_{\lambda}
\frac {\theta(\frac {x_l}{2\pi\ii}+\lambda(g)-\lambda(h)\tau-z)}
{\theta(\frac {x_l}{2\pi\ii}+\lambda(g)-\lambda(h)\tau)}
\ee^{2\pi\ii\lambda(h)z}=
\frac
{\theta(\frac{ix_l}{2\pi\ii}-iz,\frac {i\tau-s}j)}
{\theta(\frac{ix_l}{2\pi\ii},\frac {i\tau-s}j)}.
$$
\end{lemma}

\begin{proof}
First, we observe, that the set of pairs 
$(\lambda(g),\lambda(h))$ can be taken to
be the set of pairs $(\frac m{ij},\frac nj)$ such that 
$0\leq n\leq j-1$,  $0\leq m\leq ij-1$ and $m=ns\mod j$.

Let us check the transformation properties of the left hand side of 
the equation under $z\to z+1/i$. The exponential factors contribute
$$\exp(2\pi\ii\frac 1i\sum_{\lambda}\lambda(h))=
\exp(2\pi\ii\sum_{n=0}^{j-1}\frac nj)=
(-1)^{j-1}.$$
For each $n=0,\ldots,j-1$, the set of $\lambda(g)$ is given by
the fractional parts of $\frac {ns}{ij}+ \frac ki,k=0,\ldots,i-1$.
There will be exactly one such fractional part which is less than 
$\frac 1i$. The transformation $z\to z+1/i$ switches these fractions
around except for the extra $1$ for the fraction with $\lambda(g)<\frac 1i$.
As a result, we get the extra factor $(-1)^j$ from the numerator, so
overall the left hand side of the equation changes sign under 
$z\to z+\frac 1i$, as does the right hand side.

Now, let us check the transformation properties of the left hand side 
under $z\to z+\frac {i\tau-s}{ij}$. This variable change amounts to
$n\to n+1$, $m\to m+s$, which moves around the $\theta$s in the numerator,
except for the cases when new values of $m$ and $n$ fall out of their
prescribed ranges. In the case of $m$ falling out of its range, the extra 
factor required to put it back in is $(-1)$. It is easy to calculate 
the number of such occurrences, because the sum of all $m$ is going 
to change by ${ij}s$ which require $s$ switches to put into the 
correct range. So the extra factor from the switches of $m$ is
$(-1)^{s}$. In the case of $n$, it falls out of the range when 
it goes from $(j-1)$ to $j$. In this case we get $m=0\mod j$, so
$\lambda(g)=\frac ki, k=0,\ldots, i-1$.
The extra factors come from the transformation properties of $\theta$
and equal 
$$
(-1)^i\ee^{\sum_{k=0}^{i-1}(2\pi\ii(\frac {x_l}{2\pi\ii}+\frac ki-z)-
\pi\ii\tau)}=\ee^{-\pi\ii +ix_l-2\pi\ii iz-\pi\ii i\tau}.
$$
The exponential factors contribute $\exp(\pi\ii (j-1)\frac{(i\tau-s)}j)$,
so the overall factor is 
$$
\ee^{-\pi\ii +ix_l-2\pi\ii iz-\pi\ii i\tau+\pi\ii (j-1)\frac{(i\tau-s)}j
+\pi\ii s}
=
-\ee^{ix_l-2\pi\ii iz
-\pi\ii\frac{(i\tau-s)}j}
$$
which is exactly the effect of the transformation $z\to z+\frac{(is-\tau)}{ij}$
to the right hand side of the equation. 

It is straightforward to check that both sides have no poles and the same
zeroes as functions of $z$, therefore their ratio is a holomorphic elliptic 
function, hence a constant. It remains to observe that both sides equal 
$1$ for $z=0$.
\end{proof}

\begin{lemma}\label{producttwo}
$$
\prod_{\lambda\neq 0}
\frac
{\theta(\frac{D_c}{2\pi\ii}+\lambda(g)-\lambda(h)\tau)
\theta(\lambda(g)-\lambda(h)\tau-z)}
{\theta(\frac{D_c}{2\pi\ii}+\lambda(g)-\lambda(h)\tau-z)
\theta(\lambda(g)-\lambda(h)\tau)}
$$
$$=
\frac 
{\theta(\frac{iD_c}{2\pi\ii},\frac {i\tau-s}j)
\theta(\frac {D_c}{2\pi\ii}-z)
\theta(-iz,\frac {i\tau-s}j)\theta'(0)
}
{\theta(\frac{iD_c}{2\pi\ii}-iz,\frac {i\tau-s}j)
\theta(\frac {D_c}{2\pi\ii})
i\theta'(0,\frac{i\tau-s}j)\theta(-z)
}
.
$$
\end{lemma}

\begin{proof}
We use the result of Lemma \ref{product} with $x_l$ replaced by 
$D_c$ and the limit of the same calculation as $x_l\to 0$.
\end{proof}

\emph{Proof of Theorem \ref{DMVVpairs} continues.} By Lemmas \ref{product}
and \ref{producttwo}, we can rewrite $F_{i,j,s}$ as
$$
F_{i,j,s}=\sum_{\{\lambda_1,\ldots,\lambda_k\}}
\int_X \prod_l \Big(x_l
\frac{\theta(\frac{ix_l}{2\pi\ii}-iz,\frac {i\tau-s}j)}
{\theta(\frac{ix_l}{2\pi\ii},\frac {i\tau-s}j)}
\Big)
\prod_{c=1}^k\ee^{2\pi\ii ij(\alpha_c+1)\lambda_c(h)z}
$$
$$\times
\prod_{c=1}^k
\frac 
{\theta(\frac{iD_c}{2\pi\ii},\frac {i\tau-s}j)
\theta(-iz,\frac {i\tau-s}j)\theta'(0)
\theta(\frac {D_c}{2\pi\ii}+\lambda_c(g)-\lambda_c(h)\tau-ij(\alpha_c+1)z)
}
{\theta(\frac{iD_c}{2\pi\ii}-iz,\frac {i\tau-s}j)
i\theta'(0,\frac{i\tau-s}j)
\theta(\frac {D_c}{2\pi\ii}+\lambda_c(g)-\lambda_c(h)\tau)
\theta(-ij(\alpha_c+1)z)}.
$$

We will use the following lemma.
\begin{lemma}\label{ijsum}
$$
\sum_{\lambda}
\frac 
{
\theta(u+\lambda(g)-\lambda(h)\tau-v)
}
{
\theta(u+\lambda(g)-\lambda(h)\tau)
}
\ee^{2\pi\ii \lambda(h)v}
=
i
\frac 
{\theta'(0,\frac {i\tau-s}j)
\theta(-v)\theta(iu-\frac vj,\frac{i\tau-s}j)}
{\theta'(0)
\theta(-\frac vj,\frac{i\tau-s}j)
\theta(iu,\frac{i\tau-s}j)}.
$$
\end{lemma}

\begin{proof}
We use the following basic formula which is essentially contained 
in \cite{ellgeninvent}, where the right hand side converges for
$ \Im (\tau)>\Im (u)>0$.
$$
-\frac{\theta(u+z)\theta'(0)}{2\pi\ii\theta(u)\theta(z)}=
\sum_{k\in \ZZ} \frac{\ee^{2\pi\ii k u}}{1-\ee^{2\pi\ii z}
\ee^{2\pi\ii k\tau}}.
$$
We also recall the description of $(\lambda(g),\lambda(h))$
from Lemma \ref{product}. Not that the quotient depends on
the choice of $\lambda(g)\mod 1$ only, so we can assume that 
$\lambda(g)=\frac {ns}{ij}+\frac mi, m\in \ZZ/i\ZZ.$
Then,
$$
\sum_{\lambda}
\frac 
{
\theta(u+\lambda(g)-\lambda(h)\tau-v)
}
{
\theta(u+\lambda(g)-\lambda(h)\tau)
}
\ee^{2\pi\ii \lambda(h)v}
$$
$$
=\sum_{m=0}^{i-1}\sum_{n=0}^{j-1}
\frac 
{\theta(u-v+\frac mi+\frac {ns}{ij}-\frac nj \tau)}
{\theta(u+\frac mi+\frac{ns}{ij}-\frac nj \tau)}
\ee^{2\pi\ii v\frac n j}
$$
$$=
-\frac{2\pi\ii\theta(-v)}{\theta'(0)}
\sum_{m=0}^{i-1}\sum_{n=0}^{j-1}
\sum_{k\in \ZZ}
\ee^{2\pi\ii v\frac n j}
\frac{\ee^{2\pi\ii k (u+\frac mi+\frac{ns}{ij}-\frac nj \tau)}}
{1-\ee^{-2\pi\ii v}
\ee^{2\pi\ii k\tau}}
$$
$$
=-i\frac{2\pi\ii\theta(-v)}{\theta'(0)}
\sum_{n=0}^{j-1}
\sum_{k\in \ZZ}
\frac
{\ee^{2\pi\ii v\frac n j}\ee^{2\pi\ii kiu}\ee^{-2\pi\ii kn \frac{i\tau-s}j}}
{1-\ee^{-2\pi\ii v}\ee^{2\pi\ii ki\tau}}
$$
$$
=-i\frac{2\pi\ii\theta(-v)}{\theta'(0)}
\sum_{k\in \ZZ}
\frac
{
\ee^{2\pi\ii kiu}
}
{(1-\ee^{-2\pi\ii v}\ee^{2\pi\ii ki\tau})}
\frac
{(1-\ee^{2\pi\ii v}\ee^{-2\pi\ii k(i\tau-s)})}
{(1-
\ee^{2\pi\ii v\frac 1j}
\ee^{-2\pi\ii k\frac{i\tau-s}j})
}
$$
$$
=i\frac{2\pi\ii\theta(-v)}{\theta'(0)}
\sum_{k\in \ZZ}
\frac
{
\ee^{2\pi\ii kiu}\ee^{2\pi\ii v}\ee^{-2\pi\ii ki\tau}
}
{(1-
\ee^{2\pi\ii v\frac 1j}
\ee^{-2\pi\ii k\frac{i\tau-s}j})
}
$$
$$
=-i\ee^{2\pi\ii v(1-\frac 1j)}
\frac{2\pi\ii\theta(-v)}{\theta'(0)}
\sum_{k\in \ZZ}
\frac
{
\ee^{2\pi\ii k(iu-(j-1)\frac{i\tau-s}j)}
}
{(1-
\ee^{-2\pi\ii v\frac 1j}
\ee^{2\pi\ii k\frac{i\tau-s}j})
}
$$
$$
=i\ee^{2\pi\ii v(1-\frac 1j)}
\frac 
{
\theta(-v)
\theta'(0,\frac {i\tau-s}j)
\theta(iu-(j-1)\frac{i\tau-s}j-\frac vj,\frac{i\tau-s}j)
}
{\theta'(0)
\theta(iu-(j-1)\frac{i\tau-s}j,\frac{i\tau-s}j)
\theta(-\frac vj,\frac{i\tau-s}j)
}
$$
$$
=
i
\frac 
{\theta'(0,\frac {i\tau-s}j)
\theta(-v)\theta(iu-\frac vj,\frac{i\tau-s}j)}
{\theta'(0)
\theta(-\frac vj,\frac{i\tau-s}j)
\theta(iu,\frac{i\tau-s}j)}.
$$
In the above calculations the series are absolutely convergent, as long as 
$\Im(\tau)>0$ and $1>\frac {\Im (u)}{\Im(\tau)}>\frac {j-1}j$.
Then analytic continuation finishes the proof.
\end{proof}

\emph{Proof of Theorem \ref{DMVVpairs} continues.} By Lemma \ref{ijsum}
we can rewrite
$$
F_{i,j,s}=
\int_X \prod_l \Big(
\frac{x_l\theta(\frac{ix_l}{2\pi\ii}-iz,\frac {i\tau-s}j)}
{\theta(\frac{ix_l}{2\pi\ii},\frac {i\tau-s}j)}
\Big)
\prod_{c=1}^k
\frac 
{
\theta(-iz,\frac {i\tau-s}j)
\theta(\frac {iD_c}{2\pi\ii}
-i(\alpha_c+1)z,\frac{i\tau-s}j)
}
{\theta(\frac{iD_c}{2\pi\ii}-iz,\frac {i\tau-s}j)
\theta(-i(\alpha_c+1)z,\frac{i\tau-s}j).
}
$$
We notice that when we calculate $\int_X$, we only pick up the 
polynomials of degree $\dim X$ in $x_l$ and $D_c$, which allows 
us to conclude that 
$$
F_{i,j,s}(z,\tau)=Ell(X,D;iz,\frac {i\tau-s}j).
$$
We now recall that the contribution of the 
commuting pair of elements $g,h\in S_n$
to orbifold elliptic genus of $(X^n,D^{(n)})$ is
$\frac 1{n!}$ times 
the product of several $F_{i,j,s}$, each one corresponding to 
an orbit of the action of $\langle g,h\rangle$ on $\{1,\ldots,n\}$.
Every such orbit $I_m$ will have $i_m$, $j_m$ and $s_m\in\ZZ/j_m\ZZ$ 
uniquely specified.
So we have
$$
\sum_{n\geq 0}p^n {Ell}(X^n/S_n,D^{(n)}/S_n;z,\tau)=
\sum_{n\geq 0}p^n\sum_{gh=hg,g,h\in S_n} \frac 1{n!}
\prod_{I_m} F_{i_m,j_m,s_m}(z,\tau)
$$
$$=
\sum_{r\colon Z_{>0}\times Z_{>0}\to Z_{\geq 0}}
\frac {p^{\sum_{i,j}ijr(i,j)}}
{\prod_{i,j}r(i,j)!(ij)^{r(i,j)}}
\prod_{i,j} (\sum_{s=0}^{j-1}F_{i,j,s}(z,\tau)^{r(i,j)}).
$$
In this calculation we have used the fact that for $n=\sum_{i,j} ijr(i,j)$
there are 
$$\frac 1{\prod_{i,j} r(i,j)!}\frac {n!}{\prod_{i,j}((ij)!)^{r(i,j)}}
$$ 
ways to split 
$\{1,\ldots,n\}$ into groups of subsets so that there $r(i,j)$ subsets of 
``type $(i,j)$''. Then for each set of type $(i,j)$ there are 
$\frac{(ij)!}{ij}$ different ways to define the action of the $g$ and $h$
conjugate to the standard action we have discussed earlier.
We now conclude that 
$$
\sum_{n\geq 0}p^n {Ell}(X^n/S_n,D^{(n)}/S_n;z,\tau)=
\exp(\sum_{i,j>0}\sum_{s=0}^{j-1}\frac{p^{ij}}{ij}F_{i,j,s}(z,\tau))
$$
$$
=
\exp(\sum_{i,j>0}\sum_{s=0}^{j-1}
\frac{p^{ij}}{ij}Ell(X,D;iz,\frac {i\tau-s}j))
$$
$$
=\exp(\sum_{i,j>0}\sum_{m,l}\sum_{s=0}^{j-1}
c(m,l)\frac{p^{ij}}{ij}y^{il}q^{\frac {im}j}\ee^{2\pi\ii\frac {ms}j})
$$
$$
=\exp(\sum_{i,j>0}
c(mj,l)
\frac{p^{ij}}{i}y^{il}q^{im})
$$
$$
=\prod_{j=1}^\infty\prod_{m,l}\exp(c(mj,l)\sum_{i>0}\frac {p^{ij}}{i} 
y^{il}q^{im}
)
$$
$$
=\prod_{j=1}^\infty\prod_{m,l}\exp(-c(mj,l)
\ln(1-p^jy^{l}q^{m}))=\prod_{j=1}^\infty\prod_{m,l}
(1-p^jy^{l}q^{m})^{-c(mj,l)},
$$
which finishes the proof. 
\end{proof}

\begin{corollary} Let $X$ be a complex projective surface and 
$X^{(n)}$ be the Hilbert scheme of subschemes of $X$ of length $n$. 
Let $\sum_{m,l}c(m,l)y^lq^n$ 
be the elliptic genus of $X$. 
Then 
$$\sum_{n\geq 0}p^n Ell(X^{(n)};z,\tau)=
\prod_{i=1}^\infty \prod_{l,m}
\frac 1
{(1-p^iy^lq^m)^{c(mi,l)}}.$$
\label{hilbertschemes}
\end{corollary} 

\begin{proof}
By Theorem \ref{main}, the orbifold elliptic genus of the symmetric
power $X^n/S_n$ equals the elliptic genus of its crepant resolution,
which is provided by $X^{(n)}$ in the surface case.
\end{proof}

\begin{remark}
As a corollary of our work we easily deduce the analog of DMVV
conjecture for wreath products, see \cite{wang}.
\end{remark}

\section{Open questions}\label{open.section}
In this section we mention possible directions in which the results of 
this paper could be extended.

The biggest drawback of our technique is that it does not establish
elliptic genus of a Kawamata log-terminal pair as a graded dimension
of some natural vector space. In the smooth non-equivariant case
such description is provided by \eqref{elleuler}. Even more interesting
is the description of the elliptic genus as the graded dimension of the
vertex algebra which is the cohomology of the chiral de Rham complex
of \cite{MSV}, see \cite{ellgeninvent}. This is still open even 
in the non-equivariant case. This would be very interesting even in
at $q=0$ level, since it may give a vector space that realizes the 
stringy Hodge numbers of a singular variety $X$.

It would be also interesting to try to somehow extend the results of 
this paper to more general orbifolds (smooth stacks). 
The definition of orbifold elliptic genus (no divisor) 
was extended to this generality in \cite{liu}. While
our paper focuses on the global quotient case, it is possible that 
its techniques may still apply to the case of an algebraic variety with
at most quotient singularities. Indeed, the toroidal techniques 
are in some sense local. In a related remark, we do believe that the 
analog of our main theorem holds for the orbifold elliptic classes
of $(X,E,G)$ and $(X/G_1,E/G_1,G/G_1)$ where $G$ is an arbitrary 
normal subgroup of $G$.

The birational properties of elliptic genus of mean that it is 
preserved under $K$-equivalence (cf. \cite{Kawamata},\cite{clwang}).
It is conjectured in \cite{Kawamata} that $K$-equivalent varieties
have equivalent derived categories. It is therefore points to a possible
connection between elliptic classes and derived categories. It is 
however more likely that both objects are a part of a bigger structure
of a conformal field theory which somehow behaves well under 
$K$-equivalence. This is largely speculative at this point, but it 
would be interesting to define mathematically an invariant of a variety
which would encompass both its derived category and its elliptic genus.
The situation is even more murky for Kawamata log-terminal pairs, since it
is unclear what the correct definition of the derived category of the pair
may be.

A mirror symmetric analog of a resolution of singularities is a 
deformation to a smooth variety. Unfortunately, this theory is not nearly
as developed as the theory of birational morphisms. It would be interesting
to define an analog of a crepant resolution in this setting and to try
to check the invariance of the elliptic genus.

It is known that elliptic genus for smooth manifolds has a rigidity 
property. Recently, this property has been extended to the orbifold case
in \cite{liu}. It is reasonable to try to extend this property to 
the case of Kawamata log-terminal pairs. It is possible that the 
framework of pairs that consist of an orbifold and an equivariant bundle
over it, see \cite{liu}, will be useful.

It would be also interesting to see how the orbifold elliptic class of 
a singular variety $X$ compares to the Mather chern class of $X$,
see for example \cite{Fulton}. 

\section{Appendix. Assorted toric lemmas}\label{section.app}
In this appendix we collect several combinatorial statements which 
are useful in our study of toroidal morphisms.

\begin{lemma}\label{mainthetalemma}
Let $\Sigma$ be a simplicial fan in the first orthant of a lattice
$N=\oplus_{i}\ZZ e_i$. Moreover, let $\hat N$ be a sublattice of $N$
of finite index. We denote the quotient group $N/\hat N$ by $G$.
We further assume that each cone $C$ of $\Sigma$ is generated
by a part of a basis of $\hat N$. We denote by $x_i$ the linear functions
on $N_\CC$ that are dual to $e_i$. For each cone $C$ of maximum dimension
we denote by $\{x_{i;C}\}$ the linear combinations of $x_i$
which are dual to the generators of $C$. Let $a$ be a linear function
on $N$ which takes values $a_i$ on $e_i$ and values 
$a_{i;C}$ on the generators of $C$. Then
$$
\sum_{C\in \Sigma,\dim C=\rk N}
\sum_{g,h\in G}
\prod_{i}
\frac
{
\theta'(0)
\theta(\frac{x_{i;C}}{2\pi\ii}+g_{i;C}-h_{i;C}\tau-
a_{i;C})
}
{
2\pi\ii\,
\theta(\frac{x_{i;C}}{2\pi\ii}+g_{i;C}-h_{i;C}\tau)
\theta(-a_{i;C})
}
\ee^{2\pi\ii a_{i;C}h_{i;C}}
$$
$$
=
\vert N:\hat N\vert 
\prod_{i} 
\frac
{
\theta'(0)
\theta(\frac{x_{i}}{2\pi\ii}-a_i)
}
{
2\pi\ii\,
\theta(\frac{x_{i}}{2\pi\ii})
\theta(-a_i)
}
$$
where $g_{i;C}$ and $h_{i;C}$ denote rational numbers in the range 
$[0,1)$ that are fractional parts of the coordinates of 
the lifts of  $g$ and $h$ to $N$ in the basis of $C$.
\end{lemma}

\begin{proof}
We will argue by induction on $\rk N$, with $\rk N=0$ being the 
trivial base of induction.

Let's study transformation properties of both sides of the equation
under the translations $x_1\to x_1+2\pi\ii$ and $x_1\to x_1+2\pi\ii\tau$.
Under the transformation $x_1\to x_1+2\pi\ii$ the term of the sum that 
correspond to $C,g,h$ changes into the term that corresponds to 
$C,g+e_1,h$. Indeed, the coefficients of $e_1$ in the basis of the 
cone $C$ are the same as the coefficients of $x_1$ in the linear functions
$x_{i;C}$. As a result, both sides of the equation are unchanged under 
$x_1\to x_1+2\pi\ii$. Under the transformation $x_1\to x_1+2\pi\ii\tau$ the 
term that corresponds to $C,g,h$ changes into the term that corresponds to
$C,g,h-e_1$ times $\ee^{2\pi\ii a_1}$. Indeed, the extra factor comes 
from the exponential terms since $a_1$ is the difference between the  value
of $a$ on $h$ and $h-e_1$. We also observe that the terms of the product
are such that any lift of $h$ to $N$ gives the same value, so the fact that
some of the coefficients of $h-e_1$ in the basis of $C$ are not in $[0,1)$
is not a problem. Clearly, the right hand side of the equation has the
same transformation properties.

We will now show that the left hand side of the equation of the lemma 
has only simple poles at $x_1=2\pi\ii(\ZZ+\ZZ\tau)$, considered as a function
of $x_1$ with fixed generic values of other parameters. By the above 
transformation argument, it is enough to show there are no poles at 
the solutions to linear equations on $x_1$ given by $x_{i;C}=0$.
We only need to worry about such $x_{i;C}$ that define a non-coordinate 
hyperplane which corresponds to some cone $\bar C$ of dimension $\rk N-1$
in the interior of the first orthant. This cone $\bar C$ is contained in two
cones $C$ and $C'$ of maximum dimension and we argue that the contributions
of these cones to the singular part of the Laurent expansion around
$x_{i;C}=0$ cancel. Let $v_1,\ldots, v_{\rk N-1}, v$ be the generators of 
$C$ and $v_1,\ldots, v_{\rk N-1}, v'$ be the generators of $C'$
It is easy to see that $v+v'=\sum_{i=1}^{\rk N-1}c_iv_i$ 
for some integer $c_i$ and 
$$
x_{i;C}=x_{i;C'}+c_ix_{\rk N;C},~1\leq i\leq \rk N-1;~~
x_{\rk N;C}=-x_{\rk N;C'}.
$$
There are similar transformation formulas for $g_{i;C}$ and $h_{i;C}$.
The poles at $x_{\rk N;C}=0$ can be of order at most $1$, and they can
only occur in the case $g_{\rk N;C}=h_{\rk N;C}=g_{\rk N;C'}=h_{\rk N;C'}=0$.
As a result, we only need to calculate the residue at this pole. 
The residue of the term that corresponds to $C',g,h$ is equal to 
$$
-\frac 1c
\prod_{i=1}^{\rk N-1} 
\frac
{
\theta'(0)
\theta(\frac{x_{i;C}}{2\pi\ii}+g_{i;C}-h_{i;C}\tau-a_{i;C})
}
{
2\pi\ii\,
\theta(\frac{x_{i;C}}{2\pi\ii}+g_{i;C}-h_{i;C}\tau)
\theta(-a_{i;C})
}
\ee^{2\pi\ii a_{i;C}h_{i;C}}
$$
where $c$ is the coefficient of $x_1$ in $x_{\rk N; C}$ and cancels the 
residue of the term that corresponds to $C,g,h$.

By a standard argument from the theory of elliptic functions 
we conclude that the left hand side of the equation of the lemma has 
simple zeros at $x_1=2\pi\ii(a_1+\ZZ+\ZZ\tau)$ and at no other points.
Moreover, the ratio of the two sides of the equation is independent of
$x_1$. It is therefore enough to verify that the residues at $x_1=0$
of both sides are same. Only the terms with cones $C$ that have a face 
$\bar C$ of dimension $\rk N-1$ that lies in the side of the orthant 
spanned by $e_{>1}$ can contribute to the residue. We will denote the 
generator of $C$ that does not lie in $\bar C$ by $e_1$. The residue
occurs only for $g_{1;C}=h_{1;C}=0$ and then it equals 
$$
\frac 1c
\prod_{i>1}
\frac
{
\theta'(0)
\theta(\frac{x_{i;\bar C}}{2\pi\ii}+g_{i;\bar C}-h_{i;\bar C}\tau-a_{i;\bar C})
}
{
2\pi\ii\,
\theta(\frac{x_{i;\bar C}}{2\pi\ii}+g_{i;\bar C}-h_{i;\bar C}\tau)
\theta(-a_{i;\bar C})
}
\ee^{2\pi\ii a_{i;\bar C}h_{i;\bar C}}
$$
where $c$ is the coefficient of $x_1$ in $x_{1;C}$. Here we have observed
that $x_{i;C}$ restricts to $x_{i;\bar C}$ on $x_1=0$, and similarly for
$g_{i;C}$ and $h_{i;C}$. If the intersection of 
$\hat N$ and $N$  with the span of $e_{>1}$ are lattices $\hat N_1$ and 
$N_1$ respectively, then 
$$
c=\frac
{\vert N_1:\hat N_1\vert }
{\vert N:\hat N \vert }.
$$ 
It remains to apply the induction hypothesis to the fan $\Sigma_1$ induced
by $\Sigma$ on the span of $e_{>1}$.
\end{proof}

\begin{lemma}\label{firstorth}
Let $\Sigma$ be a simplicial fan in the first orthant of a lattice
$N=\oplus_{i}\ZZ e_i$. Moreover, let $\hat N$ be a sublattice of $N$
of finite index. We further assume that each cone $C$ of $\Sigma$ is generated
by a part of a basis of $\hat N$. We denote by $x_i$ the linear functions
on $N_\CC$ that are dual to $e_i$. For each cone $C$ of maximum dimension
we denote by $\{x_{i;C}\}$ the linear combinations of $x_i$
which are dual to the generators of $C$. Then
$$
\sum_{C\in \Sigma,\dim C=\rk N}\frac 1{\prod_{i=1}^{\rk N}x_{i,C}}
=\vert N:\hat N\vert \frac 1{\prod_{i=1}^{\rk N}x_{i}}.
$$
\end{lemma}

\begin{proof}
By Lemma \ref{mainthetalemma},
$$
\sum_{C\in \Sigma,\dim C=\rk N}
\sum_{g,h\in G}
\prod_{i}
\frac
{
\theta'(0)
\theta(\frac{\epsilon x_{i;C}}{2\pi\ii}+g_{i;C}-h_{i;C}\tau-
a_{i;C})
}
{
2\pi\ii\,
\theta(\frac{\epsilon x_{i;C}}{2\pi\ii}+g_{i;C}-h_{i;C}\tau)
\theta(-a_{i;C})
}
\ee^{2\pi\ii a_{i;C}h_{i;C}}
$$
$$
=
\vert N:\hat N\vert 
\prod_{i} 
\frac
{
\theta'(0)
\theta(\frac{\epsilon x_{i}}{2\pi\ii}-a_i)
}
{
2\pi\ii\,
\theta(\frac{\epsilon x_{i}}{2\pi\ii})
\theta(-a_i).
}
$$
It remains to look at the coefficient by $\epsilon^{-\rk N}$ in
the Laurent expansion of both sides around $\epsilon=0$.
\end{proof}

\begin{example}
In the case of Figure \ref{Fig1} the identity of Lemma \ref{firstorth} 
is
$$
\frac 1
{x_2(\frac{x_1-2x_2}3)}
+
\frac 1
{(\frac{2x_2-x_1}3)(\frac{2x_1-x_2}3)}
+
\frac 1
{x_1(\frac{x_2-2x_1}3)}
=\frac 3{x_1x_2}.
$$
\end{example}

\begin{lemma}\label{toricsum}
Let $\Sigma$ be a simplicial fan in a lattice $N$ such that the union of
all of its cones is a product of a subspace and a positive orthant.
In addition, we assume that all maximum-dimensional cones of $\Sigma$ 
are generated by a basis of $N$. Then 
$$
\sum_{C\in \Sigma,\dim C=\rk N} \frac 1{\prod_{i=1}^{\rk N}x_{i;C}}
=0
$$
where $x_{i;C}$ denote the basis of linear forms dual to the lattice
generators of $C$.
\end{lemma}

\begin{proof} 
By Lemma \ref{firstorth}, applied to the case $\hat N=N$, 
the function $\frac 1{\prod_{i=1}^{\rk N}x_{i;C}}$
is additive on $\Sigma$, so we can replace $\Sigma$ by any of its subdivision
with the same properties. After an appropriate subdivision,
we can assume that each cone of $\Sigma$ sits in one of the orthants 
and the support of $\Sigma$ is 
$\oplus_{i\not\in I}\RR_{\geq 0}e_i+\oplus_{i\in I}\RR e_i$ for 
some basis $\{e_i\}$  and some nonempty set $I$. Then we apply Lemma
\ref{firstorth} again to show that 
$$
\sum_{C\in \Sigma,\dim C=\rk N} \frac 1{\prod_{i=1}^{\rk N}x_{i;C}}
=\sum_{\{\sigma_i\}\in\{1,-1\}^I}\prod_{i\in I}\frac 1{\sigma_ix_i}
\prod_{i\not\in I}\frac 1{x_i}=0.
$$
\end{proof}

\end{document}